\documentclass[preprint,numbers,sort&compress]{elsarticle}
\usepackage{graphicx}
\usepackage{amsmath,amsfonts,mathtools,amsthm,mathdots}
\usepackage{pgfplots}
\usepackage{hyperref}
\pgfplotsset{compat=newest}
\usetikzlibrary{plotmarks}
\usetikzlibrary{arrows.meta}
\usepgfplotslibrary{patchplots}
\usepackage{grffile}

\newcommand{\rmi}{\mathrm{i}}
\newcommand{\bx}{\boldsymbol{x}}
\newcommand{\bu}{\boldsymbol{u}}
\newcommand{\bg}{\boldsymbol{g}}
\newcommand{\rme}{\mathrm{e}}
\newcommand{\RR}{\mathbb{R}}
\newcommand{\ZZ}{\mathbb{Z}}
\newcommand{\CC}{\mathbb{C}}
\theoremstyle{definition}
\newtheorem{example}{Example}

\DeclareMathOperator{\sech}{sech}
\date{}
\begin{document}
\begin{frontmatter}
\title{Matrix- and tensor-oriented numerical schemes
  for the evolutionary space-fractional
  complex\\
  Ginzburg--Landau
  equation}
\author[1]{Marco Caliari}
\ead{marco.caliari@univr.it}

\author[1,2]{Fabio Cassini\corref{cor1}}
\ead{fabio.cassini@univr.it, cassini@altamatematica.it}

\cortext[cor1]{Corresponding author}

\affiliation[1]{organization={Department of Computer Science,
    University of Verona},addressline={Strada Le Grazie, 15},
  postcode={37134},
  city={Verona},
country={Italy}}

\affiliation[2]{organization={Istituto Nazionale di Alta Matematica},
  addressline={Piazzale Aldo Moro, 5},
  postcode={00185},
  city={Roma},
  country={Italy}
}
\myfooter[L]{}

  \begin{abstract}
  In this manuscript, we propose matrix- and tensor-oriented
  methods for the numerical solution of the multidimensional
  evolutionary space-fractional complex Ginzburg--Landau equation.
  After a suitable spatial semidiscretization, the resulting system
  of ordinary differential equations is time integrated with stiff-resistant
  schemes. The needed actions of special matrix functions (e.g., inverse, exponential,
  and the so-called $\varphi$-functions) are efficiently computed in a direct
  way by exploiting the underlying tensor structure of the task and taking advantage
  of high performance BLAS and parallelizable pointwise operations.
  Several numerical experiments in 2D and 3D, where we
  apply the proposed technique in the context of linearly-implicit and
  exponential-type schemes, show the reliability and superiority of the
  approach against the state-of-the-art, allowing to obtain speedups which
  range from one to two orders of magnitude. Finally, we demonstrate that in our
  context a single GPU can be effectively exploited to boost the computations both
  on consumer- and professional-level hardware.
  \end{abstract}
\end{frontmatter}
  \section{Introduction}
The efficient numerical solution of evolutionary partial differential equations
(PDEs), which are at the basis of many complex physical phenomena and rarely
admit a closed-form analytical solution, is a fundamental challenge in numerical
analysis and scientific computing.
In particular, in this manuscript we focus on the following multidimensional
evolutionary space-fractional complex Ginzburg--Landau equation
(FCGLE)
  \begin{equation}\label{eq:FCGLE}
    \left\{\begin{aligned}
    \partial_t u(t,\bx)&=
    (\nu+\rmi\eta)\left(\partial_{x_1}^{\alpha_1}+\dots +\partial_{x_d}^{\alpha_d}\right)u(t,\bx)+\gamma u(t,\bx)\\
    &-(\kappa+\rmi \zeta)\lvert u(t,\bx)\rvert^2u(t,\bx)
    +s(t,\bx),&&t\in(0,T],\ \bx\in\RR^d,\\
    u(0,\bx)&=u_0(\bx),
    \end{aligned}\right.
  \end{equation}
  see, e.g., References~\cite{WH18,ZZL19,ZZS20,ZLPR20,ZZ21,ZOG21,LWL22,ZZS21,ZZS23}.
  Here, $u(t,\bx)\in\CC$ is the complex-valued unknown,
  the parameters $\nu$ and $\kappa$ are positive numbers,
  $\eta$, $\gamma$, and $\zeta$  are real numbers, and $s(t,\bx)$ represents
  the source term of the equation. The symbol
  $\partial_{x_\mu}^{\alpha_\mu}$, with
  $1<\alpha_\mu<2$, denotes the so-called \emph{Riesz fractional derivative}
  of order $\alpha_\mu$ along the direction $\mu$. We recall that, for a
  univariate function $v\colon \RR\to\CC$ and order $1<\alpha<2$,
  the Riesz fractional derivative can be defined as
  \begin{equation*}
    \partial_x^\alpha v(x)=
    \frac{\partial^\alpha}{\partial\lvert x\rvert^\alpha}v(x)=
    -\frac{1}{2\cos(\alpha\pi/2)\Gamma(2-\alpha)}
    \frac{d^2}{dx^2}\int_{-\infty}^{\infty}
    \lvert x-\xi\rvert^{1-\alpha}v(\xi)d\xi
    \end{equation*}
  and it satisfies (see Reference~\cite{ZZS23})
  \begin{equation*}
    \partial^\alpha_x v(x)  =
    -\frac{1}{2\cos(\alpha\pi/2)}(D^\alpha_{-\infty,x} v(x) + D^\alpha_{x,\infty} v(x)).
  \end{equation*}
  Here, $D^\alpha_{-\infty,x} v(x)$ and $D^\alpha_{x,\infty} v(x)$ denote the left and
  right Riemann--Liouville derivatives, respectively, that is
  \begin{equation*}
    D^\alpha_{-\infty,x} v(x) = \frac{1}{\Gamma(2-\alpha)}\frac{d^2}{dx^2}
    \int_{-\infty}^x(x-\xi)^{1-\alpha} v(\xi)d\xi
  \end{equation*}
  and
  \begin{equation*}
    D^\alpha_{x,\infty} v(x) = \frac{1}{\Gamma(2-\alpha)}\frac{d^2}{dx^2}
    \int_x^{\infty}(\xi-x)^{1-\alpha} v(\xi)d\xi.
  \end{equation*}
  We remark also that, for $\alpha \rightarrow 2$, the Riesz fractional
  derivative reduces
  to the classical second derivative.
  For analytical results such as well-posedness and other relevant properties
  of the FCGLE, we invite the reader to consult, e.g., Reference~\cite{PG13}
  and references therein.
  As frequently considered in the literature, we are interested in solutions
  to~\eqref{eq:FCGLE}  with support in
  $\Omega=(a_1,b_1)\times\dots\times(a_d,b_d)$
  or a fast decay to zero when $\lvert \bx\rvert\to\infty$.
  Therefore, it is not restrictive to require
  the homogeneous condition
  \begin{equation*}
    u(t,\bx)=0,\quad t\in[0,T],\quad \bx \in\RR^d\setminus \Omega.
  \end{equation*}

  The FCGLE is employed to describe a   wide variety of physical
  phenomena. More in detail, it has been firstly   introduced in
  Reference~\cite{TZ05} as a natural generalization of the classical
  Ginzburg--Landau equation to describe the dynamical processes in a
  medium with fractal dispersion. Later on, it has been considered,
  for instance, to model the dynamics of systems with long-range powerwise
  interaction (see Reference~\cite{TZ06}) and to investigate   wave solutions
  in a network of Hindmarsh--Rose neural model   with long-range diffusive
  couplings (see Reference~\cite{MTBBK16}).
  It is also worth noting that different fractional generalizations of the
  Ginzburg--Landau   equation exist in the literature. We mention here, e.g.,
  References~\cite{LLZ17,M18b}, where the authors consider the
  \emph{spectral fractional Laplacian} $(-\Delta)^{\alpha/2}$,
  defined for a smooth and $2\pi$-periodic function $v\colon \RR^d\to\CC$ as
  \begin{equation*}
    (-\Delta)^{\alpha/2}v(\bx)=\sum_{\boldsymbol k\in\ZZ^d}\lvert
    \boldsymbol k\rvert^{\alpha}\hat v_{\boldsymbol k}\rme^{\rmi \boldsymbol k\cdot
    \bx},
  \end{equation*}
  being $\hat v_{\boldsymbol k}$ the Fourier coefficients of $v(\boldsymbol x)$.
  We will not consider this case in the manuscript and invite an interested reader
  to consult for instance Reference~\cite{LPGSGZMCMAK20}.

  As previously mentioned, in general there is no analytical solution to
  the model under consideration. Thus, developing reliable and efficient
  numerical schemes for equation~\eqref{eq:FCGLE} is of great importance.
  Most of the techniques proposed in the literature start by performing a
  discretization of the spatial variables by spectral methods or finite
  differences, i.e., they work in the context of the method of lines.
  This, in particular, leads to a system of ordinary differential equations
  (ODEs) in the form
  \begin{equation}\label{eq:ODEs}
    \left\{\begin{aligned}
    \bu'(t)&=K\bu (t)+\bg(t,\bu(t)),&t&\in(0,T],\\
    \bu(0)&=\bu_0,
    \end{aligned}\right.
  \end{equation}
  where $\bu(t)$ is the complex-valued unknown vector collecting the degrees
  of freedom (DOF) on the computational space grid,
  $K$ is a complex-valued
  matrix discretizing the linear differential
  operator, and $\bg(t,\bu(t))$ is a nonlinear function of $\bu(t)$.
  Although in principle it would possible to employ standard explicit methods for
  the time marching of system~\eqref{eq:ODEs} (such as classical Runge--Kutta methods,
  see, e.g., Reference~\cite{LWL22}),
  the presence of the spatial differential operator and the usually
  large number of discretization points in space leads to a
  \emph{stiff} system of ODEs (where the stiffness is concentrated in $K$).
  Therefore, it is more appropriate to choose time marching schemes
  that avoid potentially prohibitive time step
  size restrictions.
  In this sense, many approaches have been proposed by the research community.
  In the context of the FCGLE,
  we mention here, among the others,
  implicit methods (see References~\cite{ZZL19,ZLPR20,ZZ21}),
  splitting schemes
  (see References~\cite{WH18,ZOG21,RDL24}), and exponential
  integrators (see References~\cite{ZZS20,ZZS23}).
  All these classes of methods enjoy favorable stability properties in terms of
  stiffness. However, considering the computational burden, they tend to be more 
  costly than standard explicit methods, since they
  require to compute the action of certain matrix functions of $K$ such as
  the inverse  of $(I-\theta K)$ and some exponential-like
  functions of $\theta K$ (being $I$ the identity matrix and $\theta$ a
  parameter depending on the time step size and on the specific numerical
  scheme). This is particularly relevant in a multidimensional framework, where
  the cost of computing matrix functions must be taken into account to obtain
  an efficient and competitive numerical scheme. For instance, in the context
  of linearly-implicit methods for two-dimensional FCGLE,
  in Reference~\cite{ZLPR20} the authors split the arising linear system according to the
  alternating direction implicit technique.
  In Reference~\cite{ZZS21}, a different two-step second-order linearly-implicit
  method is employed and the linear systems are solved by the generalized
  minimal residual method with a block circulant preconditioner.

  A suitable semidiscretization in space of problem~\eqref{eq:FCGLE} leads
  to a system of ODEs~\eqref{eq:ODEs} where the matrix $K$ has a so-called
  Kronecker sum structure (see the following section for more details).
  The presence of this feature can be suitably leveraged to boost the computation
  of the matrix functions needed in the time marching methods, see, e.g.,
  Reference~\cite{CCEOZ22} for Schr\"odinger,
  Reference~\cite{CC24ter} for (classical) complex Ginzburg--Landau,
  and References~\cite{CC24bis,DASS20} for systems of advection--diffusion--reaction
  equations leading to Turing patterns. However, to the best of our knowledge,
  an approach of this type hasn't been fully exploited in the context of the numerical
  solution of the multidimensional FCGLE.
  The main objective of this manuscript is precisely to leverage on the
  underlying Kronecker structure of the problem and its peculiar features to
  propose an equivalent matrix- and tensor-oriented formulation
  of the time integration schemes, which
  allows reaching unmatched computational efficiency levels.
  More in detail, in Section~\ref{sec:FD} we present the details of our
  choice for the spatial discretization. We then proceed in
  Section~\ref{sec:mfc} by presenting how to efficiently compute matrix
  functions in a two-dimensional context and how to generalize
  the approach to $d$ dimensions by introducing a suitable tensor framework.
  The proposed techniques can be embedded in classical time marching
  schemes. In particular, we consider linearly-implicit,
  splitting, and exponential integrators in Section~\ref{sec:numexp2} for
  some two-dimensional numerical examples.
  The approach, being heavily based on BLAS and
  parallelizable operations, is well-suited for higher-dimensional problems
  on modern computer hardware. In fact, we numerically demonstrate
  in Section~\ref{sec:numexp3}
  that the three-dimensional FCGLE can be efficiently solved on
  both laptops and workstations using the proposed technique, with consistent speedups
  when comparing GPU and CPU results.
  Finally, in Section~\ref{sec:conc} we draw the conclusions
  and mention possible future developments.
  \section{Semidiscretization of the \texorpdfstring{$d$}{d}-dimensional problem}
  \label{sec:FD}
  The semidiscretization in space of problem~\eqref{eq:FCGLE} can be performed
  in many ways, e.g., by finite differences or spectral methods.
  Similarly to Reference~\cite{CD12},
  we choose to employ a centered finite difference approximation of the spatial
  fractional differential operator. To set up the notation, we consider the
  univariate function $v$ with support in $\Omega = (a,b)$.
  We then introduce a uniform space discretization
  constituted by $n$ inner points $x^j=a+jh$, with $j=1,\dots,n$ and
  $h=(b-a)/(n+1)$ denoting the space size, and employ the approximation
  \begin{equation*}
    \partial_x^\alpha v(x^i)\approx
    -\frac{1}{h^\alpha}\sum_{j=1}^{n}g^{(\alpha)}_{\lvert j-i\rvert}
    v(x^j).
  \end{equation*}
  Here, the coefficients of the finite difference approximation are given by
  \begin{equation*}
    g_k^{(\alpha)}=\frac{(-1)^k\Gamma(\alpha+1)}{\Gamma(\alpha/2-k+1)
    \Gamma(\alpha/2+k+1)},\quad k=0,\dots,n-1,
  \end{equation*}
  where $\Gamma$ denotes the Euler gamma function.
  Notice that for $\alpha \rightarrow 2$ this approximation reduces to the
  classical second-order approximation of the second derivative. Also,
  notice that the terms $\Gamma(\alpha/2\pm k+1)$ in the denominator can
  easily over/underflow for large values of $k$. In fact, it is numerically
  more stable to compute the coefficients $g_k^{(\alpha)}$ through
  the recurrence relation (see again Reference~\cite{CD12})
  \begin{equation*}
    \left\{\begin{aligned}
      g_0^{(\alpha)}&=\frac{\Gamma(\alpha+1)}{\Gamma(\alpha/2+1)^2},\\
      g_{k}^{(\alpha)}&=\left(1-\frac{\alpha+1}{\alpha/2+k}\right)g_{k-1}^{(\alpha)},& k&=1,\dots,n-1.
    \end{aligned}\right.
  \end{equation*}
  If the function $v$ is sufficiently smooth, we can prove that the finite
  difference scheme defined
  above
  is a second-order approximation of the Riesz fractional
  derivative of $v$ (see, e.g.,
  References~\cite{CD12,SSG16,XW19,ZLPR20}), that is
  \begin{equation}\label{eq:FD2}
\partial_x^\alpha v(x^i)=-\frac{1}{h^\alpha}\sum_{j=1}^{n}g^{(\alpha)}_{\lvert j-i\rvert}
    v(x^j)+\mathcal{O}(h^2).
  \end{equation}
  If we collect the DOF $v(x^i)$ in the vector $\boldsymbol v$
  and the coefficients $g_k^{(\alpha)}$ in the matrix
  \begin{equation}\label{eq:Dnalpha}
D_{n}^{\alpha}=-\frac{1}{h^{\alpha}}\begin{bmatrix}
      g_0^{(\alpha)} & g_1^{(\alpha)} & \dots & g_{n-1}^{(\alpha)}\\
      g_1^{(\alpha)} & g_0^{(\alpha)} &\dots & g_{n-2}^{(\alpha)}\\
      \vdots & \ddots & \ddots & \vdots\\
      g_{n-1}^{(\alpha)} & \dots & g_1^{(\alpha)} & g_{0}^{(\alpha)}
      \end{bmatrix}\in \RR^{n\times n},
  \end{equation}
  we get
    $\partial_x^\alpha v(x^i) \approx (D_{n}^{\alpha}\boldsymbol v)_i$.
  Notice that the matrix~\eqref{eq:Dnalpha} is real, symmetric, Toeplitz,
  and
  negative-definite (by Gershgorin's circle theorem,
  see Reference~\cite{CD12}).
  In particular, it can be diagonalized as
  \begin{equation}\label{eq:diag}
    D_{n}^{\alpha}=Q_{n}^{\alpha}\Lambda_{n}^{\alpha}
    (Q_{n}^{\alpha})^{\mathsf T},
  \end{equation}
  where $Q_{n}^\alpha$ is an orthogonal matrix
  and $\Lambda_{n}^\alpha=\mathrm{diag}\{\lambda_1^{(\alpha)},\dots,\lambda_{n}^{(\alpha)}\}$
  is the diagonal matrix collecting the (negative) eigenvalues.
  The superscript $\sf T$ denotes the transposition.
  Also, remark that~\eqref{eq:Dnalpha} is a \emph{dense} matrix, in accordance
  with the non-locality of the fractional differential operator.

  Starting from approximation~\eqref{eq:FD2} and employing an extrapolation
  formula, we can derive a fourth-order centered finite
  difference approximation
  of the fractional differential operator.
  In fact, assuming enough smoothness,
  we have (see Reference~\cite{XW19})
  \begin{equation}\label{eq:FD4}
    \partial_x^\alpha v(x^i)=-\frac{1}{h^\alpha}\sum_{j=1}^{n}
    \check g^{(\alpha)}_{\lvert j-i\rvert}
    v(x^j)+\mathcal{O}(h^4),
  \end{equation}
  where the coefficients are defined by
  \begin{equation*}
    \check g_k^{(\alpha)}=\left\{\begin{aligned}
    &\frac{4}{3} g_k^{(\alpha)},&& \text{$k$ odd},\\
    &\frac{4}{3} g_k^{(\alpha)}-\frac{1}{3}\frac{g_{k/2}^{(\alpha)}}{2^\alpha},&& \text{$k$ even}.
    \end{aligned}\right.
  \end{equation*}
  Remark that the corresponding discretization matrix
  (denoted again $D_n^\alpha$ for simplicity of exposition) is still real,
  symmetric, Toeplitz, negative-definite, and dense. In addition, it reduces to
  the classical pentadiagonal fourth-order centered finite difference
  matrix for
  $\alpha\rightarrow 2$.

  With the introduced notation, the semidiscretization in space of
  the multidimensional problem~\eqref{eq:FCGLE} on a tensor product Cartesian grid
  is performed as follows. We start with the
  two-dimensional problem, that is $d=2$.
  We consider the spatial domain $\Omega = (a_1,b_1) \times (a_2,b_2)$
  and the inner discretization points
  $x_1^{j_1} = a_1 + j_1 h_1$  and
  $x_2^{j_2} = a_2 + j_2 h_2$ for the first and second direction, respectively.
  Here, the indexes $j_1$ and $j_2$ run from $1$ to $n_1$ and from
  $1$ to $n_2$ and the space sizes are defined as $h_1 = (b_1-a_1)/(n_1+1)$ and
  $h_2 = (b_2-a_2)/(n_2+1)$. The total number of DOF is
  $N=n_1 n_2$ and the evaluation of the unknown on the computational grid is
  denoted by $u_j(t)$, that is the $j$-th element of the vector
  $\boldsymbol u(t)\in \CC^N$
  is $u_j(t)\approx u(t,x_1^{j_1},x_2^{j_2})$, with $j = j_1 + n_1(j_2-1)$.
  We finally end up with a system of ODEs in the form~\eqref{eq:ODEs}, where the
  discretization matrix $K\in \CC^{N\times N}$ is the \emph{Kronecker sum} given by
  \begin{equation*}
    K=A_2\oplus A_1= I_2 \otimes A_1 + A_2 \otimes I_1,
  \end{equation*}
  with
  \begin{equation*}
    A_1=(\nu+\rmi \eta)D_1
    \in\CC^{n_1\times n_1} \quad \text{and}
    \quad A_2=(\nu+\rmi \eta)D_2\in\CC^{n_2\times n_2}.
  \end{equation*}
  Here and throughout the paper $\otimes$ denotes the standard Kronecker
  product between matrices, $I_\mu$ is an identity matrix of size
  $n_\mu \times n_\mu$, and we employ the shorthand notation
  $D_\mu=D_{n_\mu}^{\alpha_\mu}$.
  The $j$-th component of the vector $\bg(t,\bu(t))$ is simply given by
  \begin{equation*}
    \gamma u_j(t)-(\kappa+\rmi\zeta)\lvert u_j(t)\rvert^2u_j(t)+s(t,x_1^{j_1},x_2^{j_2}).
  \end{equation*}
  The generalization to the $d$-dimensional case is straightforward.
  Starting with the domain $\Omega = (a_1,b_1) \times \dots \times (a_d,b_d)$,
  we consider the inner discretization points
  $x_\mu^{j_\mu} = a_\mu + j_\mu h_\mu$ for $j_\mu = 1,\dots, n_\mu$,
  $h_\mu = (b_\mu-a_\mu)/(n_\mu+1)$, and $\mu=1,\dots,d$.
  Letting
  \begin{equation}\label{eq:indexj}
    j=j_1+\sum_{\mu=2}^d n_1\cdots n_{\mu-1}(j_\mu-1)
  \end{equation}
  and $u_j(t)\approx u(t,x_1^{j_1},\dots,x_d^{j_d})$ the $j$-th element of
  the vector $\boldsymbol u(t)\in \CC^N$, with $N = n_1\dots n_d$,
  we obtain a system of ODEs in the form~\eqref{eq:ODEs} with Kronecker sum
  \begin{equation}\label{eq:Ksum}
    K=A_d\oplus \dots\oplus A_1=\sum_{\mu=1}^d I_{d}\otimes\dots
    \otimes I_{\mu+1}\otimes A_\mu\otimes I_{\mu-1}\otimes\dots\otimes I_1,
  \end{equation}
  where
  \begin{equation*}
    A_\mu=(\nu+\rmi \eta)D_\mu
    \in\CC^{n_\mu\times n_\mu}.
  \end{equation*}
  The $j$-th component of the vector $\bg(t,\bu(t))$ is
  \begin{equation*}
    \gamma u_j(t)-(\kappa+\rmi\zeta)\lvert u_j(t)\rvert^2u_j(t)+s(t,x_1^{j_1},\dots,x_d^{j_d}).
  \end{equation*}

  \section{Matrix function computation}\label{sec:mfc}
  In this section, we present how to exploit the peculiar structure of the
  discretization matrix~\eqref{eq:Ksum} to accelerate the computation of
  the matrix functions needed for the time marching schemes employed
  later on. For convenience of
  the reader, we start in Section~\ref{sec:two} with the two-dimensional
  case
  and show how to equivalently rewrite it in \emph{matrix}-oriented
  formulation. We address the $d$-dimensional generalization, and the
  corresponding \emph{tensor}-oriented formulation, in Section~\ref{sec:ddim}.
  \subsection{The two-dimensional case}\label{sec:two}
  As presented above, the discretization matrix of the system for $d=2$ is
  \begin{equation*}
    K=I_2\otimes A_1+ A_2\otimes I_1,
  \end{equation*}
  with $A_\mu=(\nu+\rmi \eta) D_\mu$.
  Notice that the matrix $K\in\CC^{N\times N}$
  \emph{cannot} be considered sparse, since it
  originates from dense matrices $A_\mu$. In fact, if we set
  for simplicity $n_1=n_2=n$, $K$ has  $N^2=n^4$ elements,
  of which $\mathcal{O}(n^3)$ different from zero.
  A way to realize
  the action $K\bu$ \emph{without} assembling the matrix $K$ itself is to exploit 
  an equivalent matrix formulation of the task.
  In fact, by using the properties of the Kronecker product (see, e.g.,
  Reference~\cite[formula (2)]{VL00}) we have
  \begin{equation}\label{eq:Ku}
K\bu = \mathrm{vec}(A_1\boldsymbol U+\boldsymbol UA_2^{\mathsf T}).
  \end{equation}
  Here, $\boldsymbol U\in\CC^{n_1\times n_2}$ is such that
  $\mathrm{vec}(\boldsymbol U)=\bu \in \CC^N$, with $N=n_1n_2$, and
  $\mathrm{vec}$ is the operator which stacks the columns of the input
  matrix into a suitable column vector (that is, component $(j_1,j_2)$
  becomes component $j_1+n_1(j_2-1)$).  

Concerning the computation of functions of the matrix $K$, we first notice that
from formula~\eqref{eq:diag} we obtain the diagonalization of the matrix $A_\mu$
as
  \begin{equation*}
    A_\mu=Q_\mu\left((\nu+\rmi\eta)\Lambda_{\mu}\right)Q_\mu^{\mathsf T},
  \end{equation*}
where we denoted $Q_\mu=Q_{n_\mu}^{\alpha_\mu}$ and $\Lambda_\mu=\Lambda_{n_\mu}^{\alpha_\mu}$.
By using the mixed-product rule of the Kronecker product we easily get
$K=Q\Lambda Q^{ \mathsf T}$,
where $Q=Q_2\otimes Q_1\in\RR^{N\times N}$ is an orthogonal matrix,
$Q^{ \mathsf T} = Q_2^{ \mathsf T}\otimes Q_1^{ \mathsf T}$, and $\Lambda$ is
the complex diagonal matrix
\begin{equation*}
  \Lambda=
 (\nu+\rmi\eta)(I_{2}\otimes\Lambda_1+
 \Lambda_2\otimes I_1).
  \end{equation*}
Hence, given an analytic function $f$ and a parameter $\theta$,
we can compute
$f(\theta K)$ as
$Qf(\theta\Lambda)Q^{\mathsf T}$.
Remark that the computation of the matrix function $f$ exploiting a
diagonalization approach is well-conditioned in our case, since $Q$ is an
orthogonal matrix and hence its condition number
(in the 2-norm) is one (see also Reference~\cite[\S~4.5]{H08}).
Although we reduced the problem to a
\emph{diagonal} matrix function to compute,
forming the \emph{large-sized} and \emph{dense}
matrices $Q$ and $Q^{\sf T}$ may have a prohibitive cost.
However, we propose to express the action of $f(\theta K)$ on
the application vector $\boldsymbol v\in\CC^N$ in the matrix form
\begin{equation}\label{eq:fKu_diag2}
f(\theta K)\boldsymbol v= \mathrm{vec}(Q_1(\boldsymbol F\circ (Q_1^{\mathsf T}\boldsymbol V Q_2))Q_2^{\mathsf T}).
\end{equation}
Here, $\boldsymbol V\in\CC^{n_1\times n_2}$,
the symbol $\circ$ denotes the Hadamard product, and
$\boldsymbol F\in\CC^{n_1\times n_2}$ is obtained
collecting the scalar evaluations $f(\theta\lambda_{j_1j_2})$, where
\begin{equation*}
  \lambda_{j_1j_2}=(\nu+\rmi\eta)(\lambda_{j_1}^{(\alpha_1)}+\lambda_{j_2}^{(\alpha_2)})
\end{equation*}
are the eigenvalues of $K$ computed as direct sum of the eigenvalues
of $A_1$ and $A_2$. Also, we employed the formula
\begin{equation}\label{eq:V2kronV1u}
  (M_2\otimes M_1)\boldsymbol w = \mathrm{vec}(M_1\boldsymbol W M_2^{\mathsf T}),
  \quad M_\mu\in\CC^{n_\mu\times n_\mu},\ \boldsymbol w\in\CC^N,
  \ \boldsymbol W\in\CC^{n_1\times n_2},
\end{equation}
see again Reference~\cite[formula (2)]{VL00}.
We remark that in formula~\eqref{eq:fKu_diag2}
no large $N\times N$ matrix has to be assembled.

In our context, the proposed technique can be effectively used both in
linearly-implicit or exponential-like methods,
in which $f(z)=1/(1-z)$ or $f(z)=\varphi_\ell(z)$
(see formula~\eqref{eq:phi}).
In these cases, the evaluation cost of $f(\theta\lambda_{j_1j_2})$ is clearly
negligible.
Therefore, in formula~\eqref{eq:fKu_diag2},
the four matrix-matrix multiplications
(which can be realized by high performance level 3 BLAS)
are the largest computational cost.
In addition,  we notice that, if $f(z)=1/(1-z)$,
the described technique for the computation of $f(\theta K)\boldsymbol v$
is equivalent to the solution of the Sylvester equation
\begin{equation*}
  \left(\tfrac{I_1}{2}-\theta A_1\right)\boldsymbol X
  +\boldsymbol X\left(\tfrac{I_2}{2}-\theta A_2\right)^{\sf T}
  =\boldsymbol V,\quad \boldsymbol X\in\CC^{n_1\times n_2}
\end{equation*}
through diagonalization (see, for instance, Reference~\cite{DASS20}).
Finally, we highlight that in the special case $f(z)=\varphi_0(z)=\rme^z$,
it is preferable to
use an alternative formula to \eqref{eq:fKu_diag2}.
Indeed, taking into account that
the matrices $I_2\otimes A_1$ and $A_2\otimes I_1$ commute,
the following equivalences
\begin{equation*}
  \begin{aligned}
  \exp(\theta K)&=
  \exp(\theta (I_2\otimes A_1+A_2\otimes I_1))\\&=
  \exp(\theta (I_2\otimes A_1))\exp(\theta(A_2\otimes I_1))
  =
  \exp(\theta A_2)\otimes \exp(\theta A_1)
\end{aligned}
  \end{equation*}
hold true. Hence,
using formula~\eqref{eq:V2kronV1u},
we compute the action  $\exp(\theta K)\boldsymbol v$
as
\begin{equation}\label{eq:Neudecker}
(\exp(\theta A_2)\otimes \exp(\theta A_1))\boldsymbol v=
\mathrm{vec}(\exp(\theta A_1)\boldsymbol V\exp(\theta A_2)^{\mathsf T}).
\end{equation}
This approach is much more convenient
than~\eqref{eq:fKu_diag2} since, after the computation
of the small-sized  matrix exponentials $\exp(\theta A_\mu)$
(easily realized through the diagonalization of $A_\mu$), it only requires
two matrix-matrix multiplications.
\subsection{Extension to the \texorpdfstring{$d$}{d}-dimensional
  case}\label{sec:ddim}
We consider now the extension to $d$ dimensions of the proposed approach.
In this case, we recall that the discretization matrix of system~\eqref{eq:ODEs} is
  \begin{equation*}
    K=A_d\oplus \dots\oplus A_1=\sum_{\mu=1}^d I_{d}\otimes\dots
    \otimes I_{\mu+1}\otimes A_\mu\otimes I_{\mu-1}\otimes\dots\otimes I_1,
  \end{equation*}
  with $A_\mu=(\nu+\rmi\eta)D_\mu$.
  The generalization of formula~\eqref{eq:Ku}
to $d$ dimensions is
\begin{equation}\label{eq:kronsumv}
  K\bu=\mathrm{vec}\left(\sum_{\mu=1}^d \boldsymbol U\times_\mu A_\mu\right).
\end{equation}
Here, $\boldsymbol U\in\CC^{n_1\times \dots \times n_d}$ is an order-$d$
tensor such that
  $\mathrm{vec}(\boldsymbol U)=\bu \in \CC^N$, with $N=n_1\cdots n_d$, and
  $\mathrm{vec}$ is the operator which stacks the columns of the input
  tensor into a suitable column vector (that is,
  component $(j_1,\dots,j_d)$  becomes component $j$ according to
  \eqref{eq:indexj}).
The  symbol $\times_\mu$ denotes the so-called $\mu$-mode product, which is
the tensor-matrix operation which
multiplies a matrix onto the $\mu$-fibers of a tensor. For instance,
for $d=2$, the 1-mode product multiplies a matrix onto the
columns (1-fibers) of an order-2 tensor (i.e., a matrix) and the 2-mode
product multiplies a matrix onto the rows (2-fibers). In other words,
$\boldsymbol U\times_1 A_1=A_1\boldsymbol U$ and
$\boldsymbol U\times_2 A_2=\left(A_2\boldsymbol U^{\sf T}\right)^{\sf T}$ and
therefore formula~\eqref{eq:kronsumv} reduces to~\eqref{eq:Ku}.
For $d=3$, the 1-mode, 2-mode, and 3-mode products multiply a matrix
onto the columns, the rows, and the tubes (3-fibers)
of an order-3 tensor. We invite interested readers
to consult, e.g., References~\cite{KB09,CCZ23bis} for more details.

By using again the diagonalization of the matrices $A_\mu$ and
the mixed-product rule, the matrix $K$ can be written as
  \begin{equation*}
    K=Q\Lambda Q^{ \mathsf T},\quad Q=Q_d\otimes\dots \otimes Q_1,
\end{equation*}
where $\Lambda$ is
the complex diagonal matrix
\begin{equation*}
  \Lambda=(\nu+\rmi\eta)\sum_{\mu=1}^d I_{d}\otimes\dots
    \otimes I_{\mu+1}\otimes \Lambda_\mu
    \otimes I_{\mu-1}\otimes\dots\otimes I_1.
  \end{equation*}
We then propose to express the action of $f(\theta K)$ on the application
vector
$\boldsymbol v\in\CC^N$ in the tensor form
\begin{equation}\label{eq:fKu_diag}
  f(\theta K)\boldsymbol v = \mathrm{vec}\left(\Big(\boldsymbol F\circ\left(
  \boldsymbol V\times_1 Q_1^{\sf T}\times_2\dots\times_d Q_d^{\sf T}\right)\Big)
  \times_1 Q_1\times_2\dots\times_d Q_d\right).
\end{equation}
Here, $\boldsymbol V\in\CC^{n_1\times\dots\times n_d}$ and
$\boldsymbol F\in\CC^{n_1\times\dots\times n_d}$ is obtained
collecting the scalar evaluations $f(\theta \lambda_{j_1\dots j_d})$,
  where
\begin{equation*}
  \lambda_{j_1\dots j_d}=(\nu+\rmi\eta)\sum_{\mu=1}^d\lambda_{j_\mu}^{(\alpha_\mu)}
\end{equation*}
are the eigenvalues of $K$ computed as the direct sum of the eigenvalues of
$A_\mu$. We employed also the generalization of
formula~\eqref{eq:V2kronV1u}
\begin{equation*}
  (M_d\otimes\dots\otimes M_1)\boldsymbol w=
  \mathrm{vec}(\boldsymbol W\times_1 M_1\times_2\dots\times_d M_d),\quad
  M_\mu\in\CC^{n_\mu\times n_\mu},\ \boldsymbol w\in\CC^N,
  \ \boldsymbol W\in\CC^{n_1\times\dots\times n_d}.
\end{equation*}
We remark that the sequence of $\mu$-mode products in the
right hand side of the formula above is
known as the \emph{Tucker operator}.
Also, notice that the main computational cost
in formula~\eqref{eq:fKu_diag} is related to the $\mu$-mode products,
each one of cost $\mathcal{O}(n^{d+1})$
(if, for simplicity, we consider $n_1=\dots=n_d=n$) and
realized in practice with a single BLAS call. We invite the reader to
consult Reference~\cite{CCZ23bis} and references therein for insights on
the $d$-dimensional implementation of the Tucker operator and related operations.
The generalization of
formula~\eqref{eq:Neudecker} to the $d$-dimensional case
for the computation of $\exp(\theta K)\boldsymbol v$  is
\begin{equation}\label{eq:Neudeckerd}
  (\exp(\theta A_d)\otimes \dots
\otimes \exp(\theta A_1))\boldsymbol v =
\mathrm{vec}\left(\boldsymbol V\times_1 \exp(\theta A_1) \times_2\dots\times_d
  \exp(\theta A_d)\right).
\end{equation}
For the remaining part of the manuscript we will use the
compact notation
\begin{equation*}
  \boldsymbol W\bigtimes_{\mu=1}^d M_\mu
  =\boldsymbol W\times_1 M_1\times_2\dots\times_d M_d
\end{equation*}
for the Tucker operator with $d>2$, and keep the notation
$M_1 \boldsymbol W M_2^{\sf T}$ for $d=2$.

Finally, we highlight that, in the context of, e.g., classical
linearly-implicit or exponential-type time marching methods,
the employment of the proposed tensor-oriented technique to
evaluate matrix functions leads to
\emph{direct} schemes. In particular, there is no need of iterative methods for
the solution of linear systems or for
the approximation of matrix functions, and
therefore the computational cost per time step is essentially constant.
Also, notice that the employment of formulas~\eqref{eq:fKu_diag}
and~\eqref{eq:Neudeckerd}
allows to compute the relevant actions of matrix
functions \emph{exactly} up to machine precision, without
the need of choosing and/or tuning input tolerances or
hyperparameters.

\section{Numerical experiments}\label{sec:numexp}

In this section, we demonstrate the superiority of the proposed
matrix and tensor approaches to numerically solve the FCGLE
in a variety of two- and three-dimensional situations
using different time marching methods and software/hardware architectures.
The examples in 2D have been taken from the literature.
However, we did not find manuscripts with numerical experiments
for the Riesz fractional complex Ginzburg--Landau equation
in the 3D case. Therefore, we directly extended the 2D ones.
The selected examples, with corresponding set of parameters, are listed in
the following.
\begin{example}\label{ex:source}
  We consider the FCGLE~\eqref{eq:FCGLE}
  with parameters
  $\nu=\eta=1$, $\gamma=3$, $\kappa=1$, and $\zeta=2$. The final time is set to
  $T=1$. For the two-dimensional
  experiments, the computational domain is $\Omega=(-1,1)^2$,
  the orders of the fractional derivatives are $(\alpha_1,\alpha_2)=(1.2,1.8)$,
  and the source term $s(t,x_1,x_2)$ and the initial condition are such that
  the exact solution is
      \begin{equation*}
        u(t,x_1,x_2)=\rme^{-\rmi t}(1-x_1^2)^4(1-x_2^2)^4,
      \end{equation*}
  see~References \cite[Example 5.1]{ZZS20},
  \cite[Example 1]{ZZS23}, and \cite[Example 2]{ZLPR20}.
  For the three-dimensional experiments, we set the computational domain
  $\Omega=(-1,1)^3$, the orders of the fractional derivatives
  $(\alpha_1,\alpha_2,\alpha_3)=(1.2,1.8,1.5)$, and the source
  term $s(t,x_1,x_2,x_3)$ and the initial condition such that the exact
  solution is
  \begin{equation*}
    u(t,x_1,x_2,x_3)=\rme^{-\rmi t}(1-x_1^2)^4(1-x_2^2)^4(1-x_3^2)^4.
  \end{equation*}
\end{example}
\begin{example}\label{ex:nosource}
  We consider the FCGLE~\eqref{eq:FCGLE}
  with parameters $\nu=\eta=1$, $\gamma=1$, and $\kappa=\zeta=1$.
  The final time
  is set to $T=1$. For the two-dimensional
  experiments, the computational domain is $\Omega=(-10,10)^2$,
  the orders of the fractional derivatives are $(\alpha_1,\alpha_2)=(1.2,1.8)$,
  and there is no source term, that is, $s(t,x_1,x_2) = 0$.
  The initial condition is given by
  \begin{equation*}
    u(0,x_1,x_2)=\sech(x_1)\sech(x_2)\rme^{\rmi(x_1+x_2)},
  \end{equation*}
  see References~\cite[Example 5.2]{ZZS20} and \cite[Example 3]{ZLPR20}.
  For the three-dimensional experiments, we set the computational
  domain $\Omega=(-10,10)^3$, the orders of the fractional derivatives
  $(\alpha_1,\alpha_2,\alpha_3)=(1.2,1.8,1.5)$, and the source term
  $s(t,x_1,x_2,x_3) = 0$.
  The initial condition is
  \begin{equation*}
    u(0,x_1,x_2,x_3)=\sech(x_1)\sech(x_2)\sech(x_3)\rme^{\rmi(x_1+x_2+x_3)}.
  \end{equation*}
\end{example}

In all cases, the equation is semidiscretized in space using centered finite
differences, according to the description in Section~\ref{sec:FD},
with equal number of grid points
in each direction (that is $n_1 = \dots = n_d = n$). We employ
the second-order approximation~\eqref{eq:FD2} in most of the experiments,
with the exception of the experiment in Section~\ref{sec:ex3_2D} in which we
use the fourth-order one~\eqref{eq:FD4}.
The absolute errors in the plots are measured in the
discrete $L^2$ norm.

Concerning the time marching, we employ a constant time step size
(denoted by $\tau$) and use different schemes depending on the specific
example under consideration.
In particular, when considering Example~\ref{ex:source}, we use the following
second-order linearized multistep BDF2 scheme
(see References~\cite{WH18b,ZLPR20})
\begin{subequations}\label{eq:lbdf2-v}
\begin{equation}
  \begin{aligned}
    \bu_{n2}&=2\bu_n-\bu_{n-1}\\
    (I-\tfrac{2\tau}{3} K)\bu_{n+1}&=\tfrac{4}{3}\bu_n-\tfrac{1}{3}\bu_{n-1}+
    \tfrac{2\tau}{3}\bg(t_{n+1},\bu_{n2}),&n&\ge1,
  \end{aligned}
\end{equation}
which originates by approximating the term $\bu_{n+1}$ appearing in the standard
implicit BDF2 method with the extrapolation $\bu_{n2}=2\bu_n-\bu_{n-1}$.
The first time step is computed by the classical backward-forward Euler method,
i.e.,
\begin{equation}
  (I-\tau K)\bu_1=\bu_0+\tau\bg(0,\bu_0).
  \end{equation}
\end{subequations}
On the other hand, when considering Example~\ref{ex:nosource}, we 
employ the well-known Strang splitting scheme (also known as split-step method),
another popular second-order scheme for the solution of the FCGLE 
(see, e.g., References~\cite{WH18,RDL24}).
In this case, the ODEs system~\eqref{eq:ODEs} is split into two subproblems,
namely
  \begin{equation*}
\left\{\begin{aligned}
    \boldsymbol v'(t)&=K\boldsymbol v(t),\\
    \boldsymbol v(0)&=\boldsymbol v_0
    \end{aligned}\right.
  \end{equation*}
  and
\begin{equation*}
\left\{\begin{aligned}
  \boldsymbol w'(t)&=\gamma\boldsymbol w(t)-(\kappa+\rmi \zeta)\lvert \boldsymbol w(t)\rvert^2\boldsymbol w(t),\\
  \boldsymbol w(0)&=\boldsymbol w_0.
  \end{aligned}\right.
\end{equation*}
The exact solution of the former is clearly
$\boldsymbol v(t)=\exp(tK)\boldsymbol v_0$, while
the exact solution of the latter is given by
  \begin{equation*}
    \boldsymbol w(t) = \rme^{t\gamma - \frac{\kappa+\rmi\zeta}{2\kappa}\log(1+2t\varphi_1(2\gamma t)\kappa\lvert \boldsymbol w_0 \rvert^2)}\boldsymbol w_0=\phi_t(\boldsymbol w_0),
  \end{equation*}
  where $\phi_t$ denotes the exact flow.
Overall, the time marching of~\eqref{eq:ODEs} with the Strang
splitting scheme is
\begin{equation}\label{eq:Strang}
\bu_{n+1}=\phi_{\tau/2}\left(\exp(\tau K)\phi_{\tau/2}(\bu_n)\right).
\end{equation}
Finally, we also consider the fourth-order exponential integrator of
Runge--Kutta type proposed by Krogstad in
Reference~\cite[formula (51)]{Kr05}, i.e.,
\begin{subequations}\label{eq:krogstad-v}
\begin{equation}
  \begin{aligned}
    \bu_{n2}&=\bu_n+\tfrac{\tau}{2}\varphi_1(\tfrac{\tau}{2}K)\boldsymbol f_n\\
    \bu_{n3}&=\bu_n+\tfrac{\tau}{2}\varphi_1(\tfrac{\tau}{2}K)\boldsymbol f_n+\tau\varphi_2(\tfrac{\tau}{2} K)\boldsymbol d_{n2}\\
    \bu_{n4}&=\bu_n+\tau\varphi_1(\tau K)\boldsymbol f_n+
    2\tau\varphi_2(\tau K)\boldsymbol d_{n3}\\
        \bu_{n+1}&=\bu_n+\tau\varphi_1(\tau K)\boldsymbol f_n+
        \tau\varphi_2(\tau K)\left(2\boldsymbol d_{n2}+2\boldsymbol d_{n3}-
        \boldsymbol d_{n4}\right)\\
        &\phantom{=}+\tau\varphi_3(\tau K)\left(-4\boldsymbol d_{n2}-4\boldsymbol d_{n3}
    +4\boldsymbol d_{n4}\right),
  \end{aligned}
\end{equation}
see also References~\cite{ZZS20,ZZS23} for its application in the context of
the FCGLE.
Here,
\begin{equation}
 \boldsymbol f_n=K\bu_n+\bg(t_n,\bu_n),\quad
  \boldsymbol d_{ni}=\bg(t_n+c_i\tau,\bu_{ni})-\bg(t_n,\bu_n),
  \ c_2=c_3=\tfrac{1}{2},\ c_4=1,
\end{equation}
\end{subequations}
and the appearing matrix functions are the so-called $\varphi$-functions,
defined for a scalar argument $z\in\CC$ by
\begin{equation}\label{eq:phi}
  \varphi_{\ell}(z)=
  \left\{
  \begin{aligned}
  &\rme^{z}, &\ell& = 0,\\
  &\frac{\varphi_{\ell-1}(z)-\frac{1}{(\ell-1)!}}{z}, &\ell&\ge 1,
  \end{aligned}
  \right.
\end{equation}
and $\varphi_{\ell}(0)=1/\ell!$.

\subsection{Two-dimensional experiments}\label{sec:numexp2}
We consider in this section the numerical experiments in two space
dimensions. We start by presenting the proposed matrix-oriented technique to 
realize the time marching methods and briefly discuss the selected terms
of comparisons.

First of all, we summarize our technique to implement the
linearized BDF2 scheme~\eqref{eq:lbdf2-v}. To this aim,
we perform the diagonalization of the relevant matrices and compute
\begin{equation*}
  \boldsymbol R_\theta=\left((1-\theta\lambda_{j_1j_2})^{-1}\right)_{j_1j_2}\in
  \CC^{n\times n},\quad \theta\in\left\{\tfrac{2\tau}{3},\tau\right\}
\end{equation*}
once and for all before the time integration starts. Then,
exploiting formula~\eqref{eq:fKu_diag2}, the numerical scheme
is implemented as
\begin{subequations}\label{eq:lbdf2-m}
\begin{equation}
    \boldsymbol U_{1}=Q_1\left(\boldsymbol R_\tau\circ \left(Q_1^{\mathsf T}\left(\boldsymbol U_0+
    \tau\boldsymbol G(0,\boldsymbol U_0)\right)Q_2\right)\right)
    Q_2^{\mathsf T}
\end{equation}
for the first time step and as
\begin{equation}
  \begin{aligned}
    \boldsymbol U_{n2}&=2\boldsymbol U_n-\boldsymbol U_{n-1},\\
  \boldsymbol U_{n+1}&=Q_1\left(\boldsymbol R_{2\tau/3}\circ \left(Q_1^{\mathsf T}\left(\tfrac{4}{3}\boldsymbol U_n-
    \tfrac{1}{3}\boldsymbol U_{n-1}+
    \tfrac{2\tau}{3}\boldsymbol G(t_{n+1},
    \boldsymbol U_{n2})\right)Q_2\right)\right)
    Q_2^{\mathsf T},
    \end{aligned}
\end{equation}
\end{subequations}
for the following ones.
We will refer to this scheme as \textsc{lbdf2-m}.
Remark that the main computational cost of the procedure at each time step is given by
the application of two 2D Tucker operators (i.e., four matrix-matrix products).
As a term of comparison, we consider the realization of the vector-oriented form of
scheme~\eqref{eq:lbdf2-v} in which the linear systems are solved
by the preconditioned generalized minimal residual (PGMRES) method.
As preconditioner, similarly to Reference~\cite{ZZS23}, in our experiments
we obtained best results employing the so-called tau-preconditioner~\cite{BB90}.
For convenience of the reader, we briefly outline here the main points.
The iterative method for the solution of the arising linear systems
requires the action of the matrix
\begin{equation*}
  I-\theta K=
  I-\theta(\nu+\rmi\eta)(I_2\otimes D_{1}+
  D_{2}\otimes I_1).
\end{equation*}
Notice that $I-\theta K$ is a block Toeplitz matrix, where the blocks are
Toeplitz matrices themselves (such matrices are usually referred to as BTTB).
Hence, we can compute its action on a vector in
$\mathcal{O}(n^2\log n)$ operations using fast Fourier transform-based techniques.
An effective preconditioner for the matrix $D_\mu$ is the aforementioned
tau-preconditioner, which is defined as
$\tau(D_{\mu})=D_{\mu}-B_{\mu}$ being $B_{\mu}$ a Hankel matrix depending on $D_{\mu}$.
Then, the tau-preconditioner for the iteration matrix is
obtained as
\begin{equation*}
  I-\theta(\nu+\rmi\eta)(I_2\otimes \tau(D_{1})+
  \tau(D_{2})\otimes I_1).
\end{equation*}
This preconditioner can be diagonalized by the sine transform, and hence
can be again applied in $\mathcal{O}(n^2\log n)$ operations.
We will refer to the realization of scheme~\eqref{eq:lbdf2-v} using this
approach as \textsc{lbdf2-v}.

Our proposed matrix-oriented technique to efficiently realize the
Strang splitting scheme~\eqref{eq:Strang} is given by
\begin{equation}\label{eq:strang-m}
\boldsymbol U_{n+1}=\Phi_{\tau/2}(E_1\Phi_{\tau/2}(\boldsymbol U_n)E_2^{\mathsf T}).
\end{equation}
Here, we employed the matrices $E_\mu=\exp(\tau A_\mu)$ (which can be
computed once and for all at the beginning since the time step size is constant)
and we used the notation $\mathrm{vec}(\Phi_{\tau/2}(\boldsymbol U_n))=\phi_{\tau/2}(\bu_n)$.
We will refer to this scheme as \textsc{strang-m}.
Remark that the main computational cost per time step is given by the application
of one 2D Tucker operator (i.e., two matrix-matrix products).
As term of comparison, we consider the vector-oriented scheme~\eqref{eq:Strang}
in which the needed action of the matrix exponential is realized by a
shift-and-invert Lanczos approach. This technique has been employed in
Reference~\cite{ZZS23} for generic $\varphi_\ell$ functions in the context of
exponential Runge--Kutta methods for the two-dimensional FCGLE,
and we briefly recall it in the following.
The classical Lanczos algorithm is used to produce
the approximation
\begin{equation*}
  \left(I-\xi(I_2\otimes D_{1}+
  D_{2}\otimes I_1)\right)^{-1}V_m\approx V_m T_m,
\end{equation*}
where $m \ll N$ is the size of the Krylov subspace, $V_m\in\CC^{N\times m}$ has
orthonormal columns, $T_m\in\RR^{m\times m}$
is a tridiagonal symmetric positive definite matrix,
and $\xi>0$ is the shift-and-invert hyperparameter. Each of the $m$ iterations
of the procedure requires solving a linear system with the
\emph{real} symmetric positive definite matrix $I-\xi(I_2\otimes D_{1}+
D_{2}\otimes I_1)$ and \emph{complex} right hand side. This is realized in
practice using the preconditioned conjugate gradient (PCG) method
(with tau-preconditioner, similarly as above). Finally, we retrieve
$\varphi_\ell(\theta K)\boldsymbol v$ as
\begin{equation*}
  \varphi_\ell(\theta K)\boldsymbol v\approx
  \lVert\boldsymbol v\rVert_2V_m\varphi_\ell\left(\tfrac{\theta}{\xi}(\nu+\rmi\eta)(I-T_m^{-1})\right)\boldsymbol e_1,
\end{equation*}
where $\boldsymbol e_1$ is the first canonical vector.
The small-sized matrix function on the right
hand side is computed through the diagonalization of $T_m$. We will
refer to as \textsc{strang-v} the embedding of this approach
into~\eqref{eq:Strang}.

Finally, for our proposed matrix-oriented implementation of the Krogstad
scheme~\eqref{eq:krogstad-v}, we perform the relevant diagonalizations and
compute the matrices
  \begin{equation*}
  \boldsymbol P_{\ell,\theta}=\left(\varphi_\ell\left(\theta
  \lambda_{j_1j_2}\right)\right)_{j_1j_2},\quad \ell\in\{1,2,3\},
  \ \theta\in\left\{\tfrac{\tau}{2},\tau\right\}
  \end{equation*}
once and for all before the actual time integration (since $\tau$ is constant).
The scheme is then realized as
\begin{subequations}\label{eq:krogstad-m}
    \begin{equation}
    \begin{aligned}
      \boldsymbol U_{n2}&=\boldsymbol U_n+
      \tfrac{\tau}{2}Q_1\left(\boldsymbol P_{1,\tau/2}
      \circ \hat {\boldsymbol F}_n\right)Q_2^{\sf T}\\
      \boldsymbol U_{n3}&=\boldsymbol U_n+
      \tau Q_1\left(\tfrac{1}{2}\boldsymbol P_{1,\tau/2}
      \circ \hat {\boldsymbol F}_n+
      \boldsymbol P_{2,\tau/2}\circ
      \hat {\boldsymbol D}_{n2}\right)Q_2^{\sf T}\\
      \boldsymbol U_{n4}&=\boldsymbol U_n+
      \tau Q_1\left(\boldsymbol P_{1,\tau}
      \circ \hat {\boldsymbol F}_n+2
      \boldsymbol P_{2,\tau}\circ
      \hat {\boldsymbol D}_{n3}\right)Q_2^{\sf T}\\
      \boldsymbol U_{n+1}&=\boldsymbol U_n+
      \tau Q_1\Big(\boldsymbol P_{1,\tau}\circ
      \hat {\boldsymbol F}_n+
      \boldsymbol P_{2,\tau}\circ \left(
      2\hat {\boldsymbol D}_{n2}+2\hat{\boldsymbol D}_{n3}-
      \hat{\boldsymbol D}_{n4}\right)\\
      &\phantom{=}+\boldsymbol P_{3,\tau}\circ \left(
      -4\hat {\boldsymbol D}_{n2}-4\hat{\boldsymbol D}_{n3}+
      4\hat{\boldsymbol D}_{n4}\right)
      \Big)Q_2^{\sf T},
    \end{aligned}
  \end{equation}
where
  \begin{equation}
    \hat {\boldsymbol F}_{n}=
Q_1^{\sf T}\left(A_1\boldsymbol U_n+\boldsymbol U_nA_2^{\sf T}+
    \boldsymbol G(t_n,\boldsymbol U_n)\right)Q_2
  \end{equation}
and
  \begin{equation}
    \hat {\boldsymbol D}_{ni}=Q_1^{\sf T}\left(
    \boldsymbol G(t_n+c_i\tau,\boldsymbol U_{ni})-
    \boldsymbol G(t_n,\boldsymbol U_n)\right)Q_2.
  \end{equation}
  \end{subequations}
We will call this approach \textsc{krogstad-m}. Remark that, at each time step,
we compute the quantities $\hat{\boldsymbol F}_n$
(one linear operator evaluation and one
2D Tucker operator) and $\hat{\boldsymbol D}_{ni}$
(one 2D Tucker operator each). Also, for the 
evaluation of each intermediate stage and the final solution we
apply an additional 2D Tucker operator.
As a term of comparison, we implement the vector-oriented version of
scheme~\eqref{eq:krogstad-v} where the actions of $\varphi$-functions are
realized by means of the shift-and-invert Lanczos technique described above
(as in Reference~\cite{ZZS23}). We will label this approach \textsc{krogstad-v}.

To perform in practice the numerical experiments, we employ
MathWorks MATLAB\textsuperscript{\textregistered} R2024b as software,
while the hardware is a standard laptop equipped with an Intel
Core i7-10750H CPU (6 physical cores) and 16GB of RAM.
The source code to reproduce the numerical experiments (fully compatible
with GNU Octave) is publicly available in a GitHub
repository\footnote{Available at \url{https://github.com/cassinif/fcgle_mt} commit \texttt{57a673a}}.
We finally highlight that the matrix diagonalizations needed for our proposed
techniques are performed through the internal MATLAB command \verb+eig+, while
for the vector-oriented approaches we exploited the routines \verb+gmres+,
\verb+pcg+, and \verb+fft2+.

  \subsubsection{First numerical experiment}\label{sec:ex1_2D}
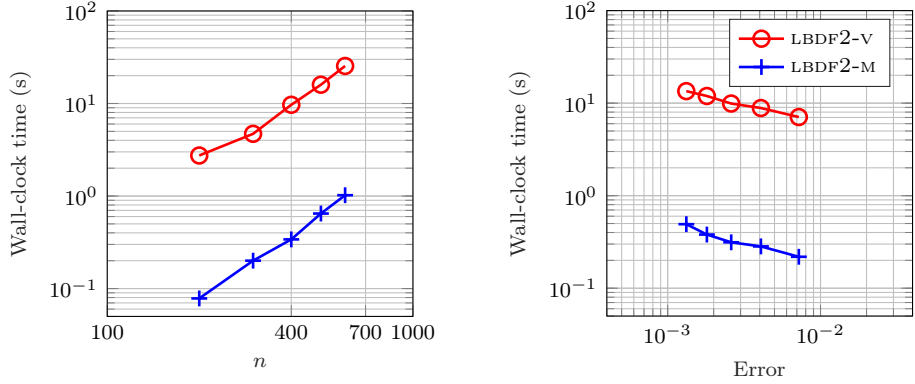
\begin{figure}[!htb]
    \centering
    \begin{tikzpicture}

\begin{axis}[%
width=1.6in,
height=1.6in,
font = \footnotesize,
scale only axis,
xmode=log,
xmin=100,
xmax=1000,
xminorticks=true,
xlabel style={font=\color{white!15!black}},
xtick={100,400,700,1000},
xticklabels={100,400,700,1000},
ymode=log,
ymin=5e-2,
ymax=1e2,
yminorticks=true,
ylabel style={font=\color{white!15!black}},
xlabel = {\footnotesize $n$},
ylabel = {\footnotesize Wall-clock time (s)},
xminorgrids,
xmajorgrids,
yminorgrids,
ymajorgrids,
axis background/.style={fill=white},
legend style={at={(0.45,0.75)}, anchor=south west, legend cell align=left, align=left, draw=white!15!black, font=\footnotesize},
]
\addplot [color=red, line width=1pt, mark size = 3pt,mark=o, mark options={solid, red}]
  table[row sep=crcr]{%
200 	2.7434e+00 \\
300 	4.7024e+00 \\
400 	9.6871e+00 \\
500 	1.5966e+01 \\
600 	2.5444e+01 \\
};

\addplot [color=blue, line width=1pt, mark size = 3pt, mark=+, mark options={solid, blue}]
  table[row sep=crcr]{%
200 	 7.8283e-02\\
300 	 2.0028e-01\\
400 	 3.3993e-01\\
500 	 6.4704e-01\\
600 	 1.0207e+00\\
};

\end{axis}

\end{tikzpicture}%
\hfill
    \raisebox{-0.8ex}{\begin{tikzpicture}

\begin{axis}[%
width=1.6in,
height=1.6in,
font = \footnotesize,
scale only axis,
xmin=4e-4,
xmax=4e-2,
xmode=log,
xminorticks=true,
xlabel style={font=\color{white!15!black}},
ymode=log,
ymin=5e-2,
ymax=1e2,
yminorticks=true,
ylabel style={font=\color{white!15!black}},
xlabel = {\footnotesize Error},
ylabel = {\footnotesize Wall-clock time (s)},
xminorgrids,
xmajorgrids,
yminorgrids,
ymajorgrids,
axis background/.style={fill=white},
legend style={at={(0.4,0.75)}, anchor=south west, legend cell align=left, align=left, draw=white!15!black, font=\footnotesize},
]
\addplot [color=red, line width=1pt, mark size = 3pt,mark=o, mark options={solid, red}]
  table[row sep=crcr]{%
7.2132e-03 	 7.0762e+00\\
4.0747e-03 	 8.8548e+00\\
2.6044e-03 	 9.9006e+00\\
1.8053e-03 	 1.1944e+01\\
1.3243e-03 	 1.3479e+01\\
};
\addlegendentry{\textsc{lbdf2-v}}

\addplot [color=blue, line width=1pt, mark size = 3pt, mark=+, mark options={solid, blue}]
  table[row sep=crcr]{%
7.2132e-03 	 2.1859e-01\\
4.0747e-03 	 2.8298e-01\\
2.6049e-03 	 3.1342e-01\\
1.8058e-03 	 3.8001e-01\\
1.3247e-03 	 4.9193e-01\\
};
\addlegendentry{\textsc{lbdf2-m}}

\end{axis}

\end{tikzpicture}%
}
    \caption{Results of Example~\ref{ex:source} in 2D using \textsc{lbdf2-m}
    and \textsc{lbdf2-v},
    see Section~\ref{sec:ex1_2D}. Left plot: increasing number of DOF per
    direction $n=[200,300,400,500,600]$, fixed number of time steps
    25. Right plot: increasing number of time steps $[15,20,25,30,35]$, fixed
    number of DOF per direction $n=400$.}
    \label{fig:ex2D1}
  \end{figure}

  We consider here Example~\ref{ex:source}, and use as time integrator
  the linearized BDF2 scheme with realization in matrix form
  (\textsc{lbdf2-m}~\eqref{eq:lbdf2-m})
  and in vector form (\textsc{lbdf2-v}), see the discussion above. Notice
  that, for \textsc{lbdf2-v}, we use in \verb+gmres+ the default tolerance
  (that is $10^{-6}$), we set the maximum number of iterations to 20,
  and we choose
  the numerical
  approximation $\bu_n$
  as initial guess for the solution vector $\bu_{n+1}$.
  The results of the experiment,
  for increasing number of DOF and increasing number of time steps,
  are summarized in Figure~\ref{fig:ex2D1}.
  In the left plot, the average
  number of PGMRES iterations for \textsc{lbdf2-v} turns out to be about 3,
  while for the right one it spans from 4 to 3.
  As we can see, the proposed technique outperforms the vector-oriented one in
  all instances, with an average speedup in terms of computational time of
  about 30. 

  \subsubsection{Second numerical experiment}\label{sec:ex2_2D}
\begin{figure}[!htb]
    \centering
    \begin{tikzpicture}

\begin{axis}[%
width=1.6in,
height=1.6in,
font = \footnotesize,
scale only axis,
xmode=log,
xmin=200,
xmax=1700,
xminorticks=true,
xlabel style={font=\color{white!15!black}},
xtick={200,700,1200,1700},
xticklabels={200,700,1200,1700},
ymode=log,
ymin=1e-1,
ymax=1e3,
yminorticks=true,
ylabel style={font=\color{white!15!black}},
xlabel = {\footnotesize $n$},
ylabel = {\footnotesize Wall-clock time (s)},
xminorgrids,
xmajorgrids,
yminorgrids,
ymajorgrids,
axis background/.style={fill=white},
legend style={at={(0.45,0.75)}, anchor=south west, legend cell align=left, align=left, draw=white!15!black, font=\footnotesize},
]
\addplot [color=orange, line width=1pt, mark size = 3pt,mark=x, mark options={solid, orange}]
  table[row sep=crcr]{%
400 	 2.2776e+01\\
600 	 6.7388e+01\\
800 	 1.4036e+02\\
1000 	 2.0275e+02\\
1200 	 3.6398e+02\\
};

\addplot [color=purple, line width=1pt, mark size = 3pt, mark=asterisk, mark options={solid, purple}]
  table[row sep=crcr]{%
400 	 1.9649e-01\\
600 	 4.3127e-01\\
800 	 1.0181e+00\\
1000 	 1.8310e+00\\
1200 	 2.8444e+00\\
};

\end{axis}

\end{tikzpicture}%
\hfill
    \raisebox{-0.8ex}{\begin{tikzpicture}

\begin{axis}[%
width=1.6in,
height=1.6in,
font = \footnotesize,
scale only axis,
xmin=4e-4,
xmax=4e-2,
xmode=log,
xminorticks=true,
xlabel style={font=\color{white!15!black}},
ymode=log,
ymin=1e-1,
ymax=1e3,
yminorticks=true,
ylabel style={font=\color{white!15!black}},
xlabel = {\footnotesize Error},
ylabel = {\footnotesize Wall-clock time (s)},
xminorgrids,
xmajorgrids,
yminorgrids,
ymajorgrids,
axis background/.style={fill=white},
legend style={at={(0.4,0.78)}, anchor=south west, legend cell align=left, align=left, draw=white!15!black, font=\footnotesize},
]
\addplot [color=orange, line width=1pt, mark size = 3pt,mark=x, mark options={solid, orange}]
  table[row sep=crcr]{%
 1.7523e-02	5.2370e+01 \\
 4.5158e-03	9.6804e+01 \\
 2.0055e-03	1.3966e+02 \\
 1.1161e-03	1.7280e+02 \\
 7.0741e-04	2.2561e+02 \\
};
\addlegendentry{\textsc{strang-v}}

\addplot [color=purple, line width=1pt, mark size = 3pt, mark=asterisk, mark options={solid, purple}]
  table[row sep=crcr]{%
1.7520e-02 	 4.3512e-01\\
4.5181e-03 	 7.4133e-01\\
2.0154e-03 	 1.0098e+00\\
1.1314e-03 	 1.2944e+00\\
7.2082e-04 	 1.5627e+00\\
};
\addlegendentry{\textsc{strang-m}}

\end{axis}

\end{tikzpicture}%
}
    \caption{Results of Example~\ref{ex:nosource} in 2D using \textsc{strang-m}
    and \textsc{strang-v},
    see Section~\ref{sec:ex2_2D}. Left plot: increasing number of DOF per
    direction $n=[400,600,800,1000,1200]$, fixed number of time steps
    15. Right plot: increasing number of time steps $[5,10,15,20,25]$, fixed
    number of DOF per direction $n=800$.}
    \label{fig:ex2D2}
  \end{figure}
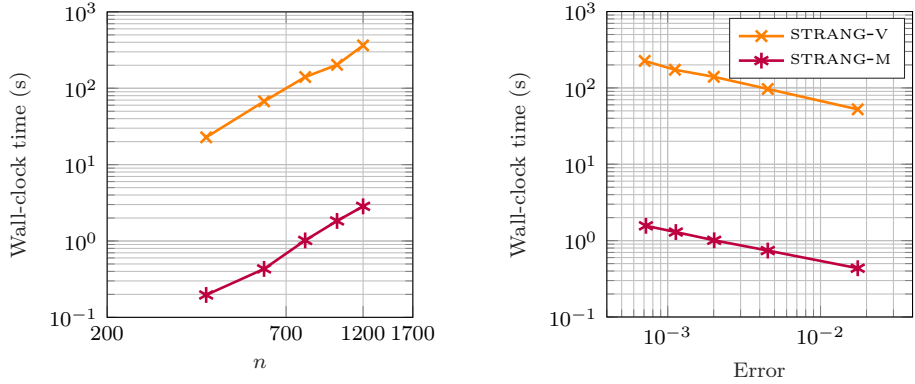

\begin{figure}[!htb]
  \centering
  \includegraphics[scale=0.475]{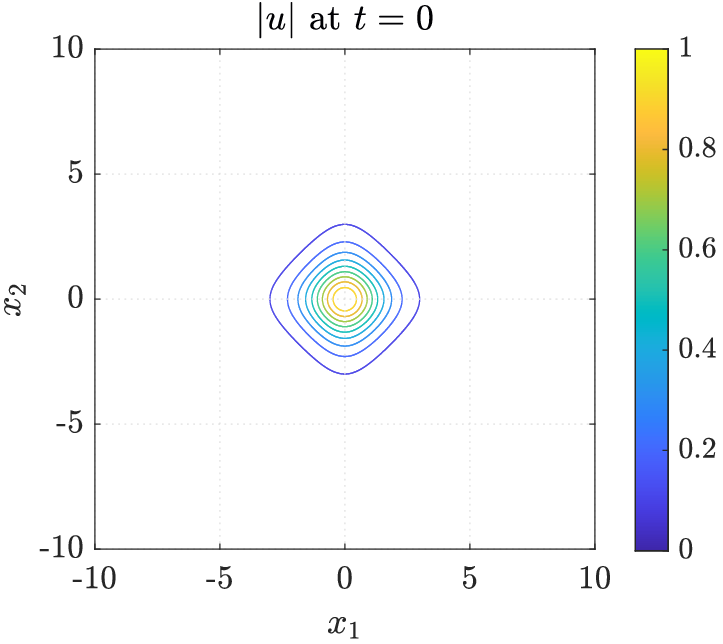}\hfill
  \includegraphics[scale=0.475]{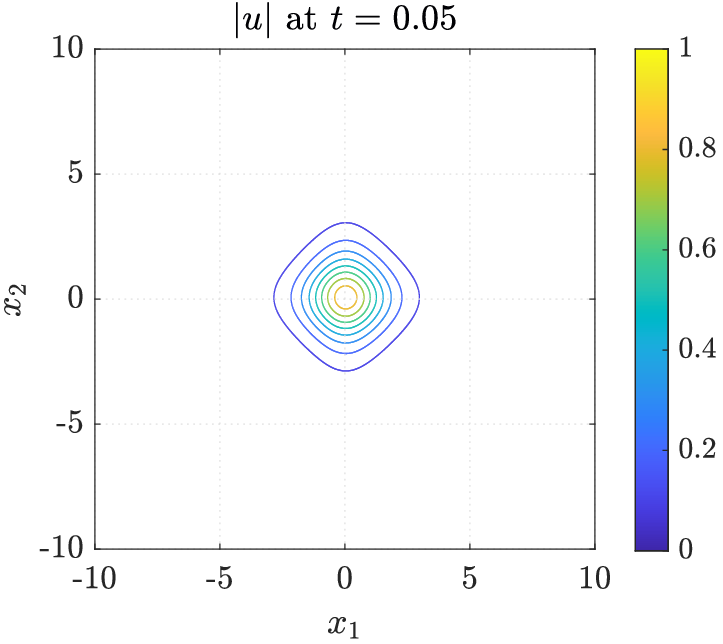}\\[2ex]
  \includegraphics[scale=0.475]{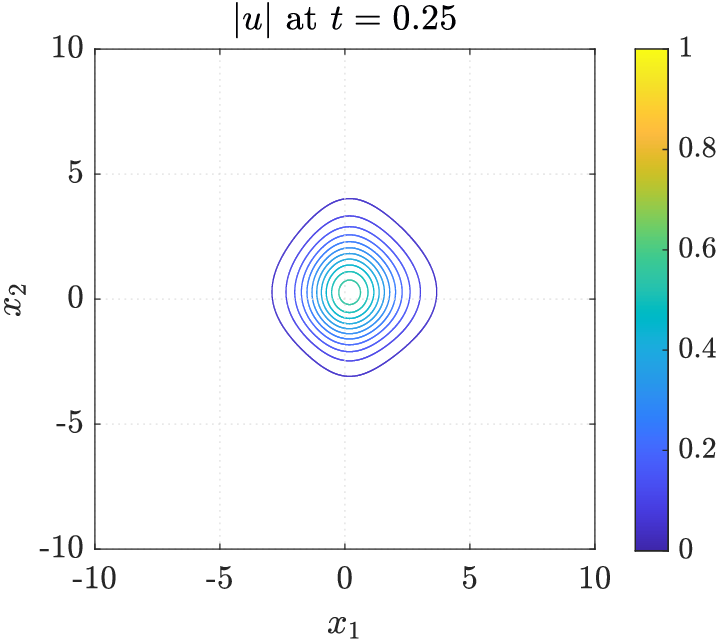}\hfill
  \includegraphics[scale=0.475]{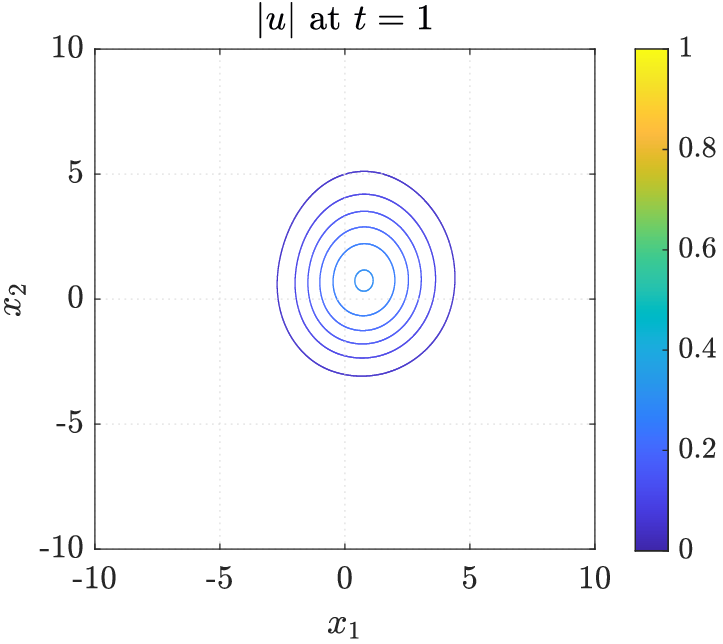}
  \caption{Level curves of $\lvert u\rvert$ at different times
  for Example~\ref{ex:nosource}
  in 2D with time integrator \textsc{strang-m}, see Section~\ref{sec:ex2_2D}.}
  \label{fig:contour_2D}
\end{figure}

  We now proceed with Example~\ref{ex:nosource}, employing as time integrator
  the Strang splitting scheme with realization in matrix form
  (\textsc{strang-m}~\eqref{eq:strang-m})
  and in vector form (\textsc{strang-v}). Remark that, for the latter, we choose
  the hyperparameters for shift-and-invert Lanczos as $m=10$ and $\xi=\tau/10$,
  and \verb+pcg+ is run with default tolerance (that is $10^{-6}$),
  maximum number of iterations 20, and initial guess the null vector.
  The results of the experiment,
  for increasing number of DOF and increasing number of time steps,
  are summarized in Figure~\ref{fig:ex2D2}.
  For the vector-oriented method, remark   that the number of PCG
  iterations spans (in average) from $2.78$ to $3.30$ for the left plot
  and from $3.58$ to $2.77$ for the right one.
  Also in this case, we demonstrate the superiority of the proposed
  matrix-oriented technique
  in all instances. As a matter of fact, the average speedup
  is more than two orders of magnitude.
  Since for Example~\ref{ex:nosource} the analytical solution is not available,
  we display in Figure~\ref{fig:contour_2D} the evolution of
  the numerical solution
  obtained with \textsc{strang-m} at different times. The results are perfectly
  in line with those available in the literature, see, e.g.,
  Reference~\cite[Fig.~5]{ZLPR20}.

\subsubsection{Third numerical experiment}\label{sec:ex3_2D}
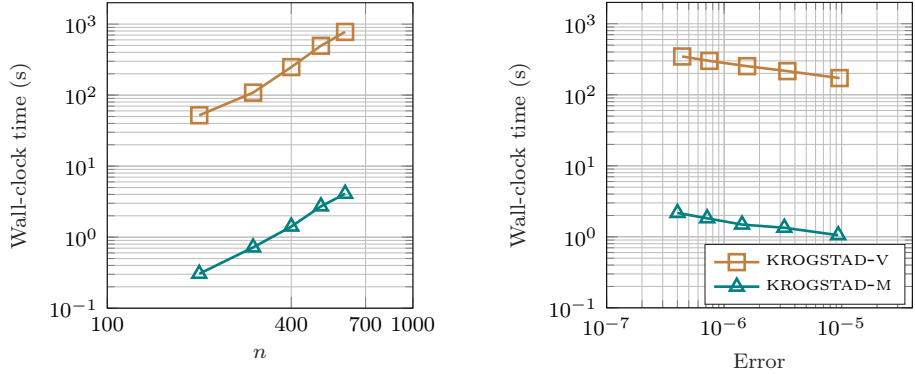
\begin{figure}[!htb]
    \centering
    \begin{tikzpicture}

\begin{axis}[%
width=1.6in,
height=1.6in,
font = \footnotesize,
scale only axis,
xmode=log,
xmin=100,
xmax=1000,
xminorticks=true,
xlabel style={font=\color{white!15!black}},
xtick={100,400,700,1000},
xticklabels={100,400,700,1000},
ymode=log,
ymin=1e-1,
ymax=2e3,
yminorticks=true,
ylabel style={font=\color{white!15!black}},
xlabel = {\footnotesize $n$},
ylabel = {\footnotesize Wall-clock time (s)},
xminorgrids,
xmajorgrids,
yminorgrids,
ymajorgrids,
axis background/.style={fill=white},
legend style={at={(0.45,0.75)}, anchor=south west, legend cell align=left, align=left, draw=white!15!black, font=\footnotesize},
]
\addplot [color=brown, line width=1pt, mark size = 3pt,mark=square, mark options={solid, brown}]
  table[row sep=crcr]{%
200 	 5.2034e+01\\
300 	 1.0846e+02\\
400 	 2.4816e+02\\
500 	 4.9509e+02\\
600 	 7.7656e+02\\
};

\addplot [color=teal, line width=1pt, mark size = 3pt, mark=triangle, mark options={solid, teal}]
  table[row sep=crcr]{%
200 	 3.0511e-01\\
300 	 7.1727e-01\\
400 	 1.4114e+00\\
500 	 2.7133e+00\\
600 	 4.0843e+00\\
};

\end{axis}

\end{tikzpicture}%
\hfill
    \raisebox{-0.8ex}{\begin{tikzpicture}

\begin{axis}[%
width=1.6in,
height=1.6in,
font = \footnotesize,
scale only axis,
xmin=1e-7,
xmax=4e-5,
xmode=log,
xminorticks=true,
xlabel style={font=\color{white!15!black}},
ymode=log,
ymin=1e-1,
ymax=2e3,
yminorticks=true,
ylabel style={font=\color{white!15!black}},
xlabel = {\footnotesize Error},
ylabel = {\footnotesize Wall-clock time (s)},
xminorgrids,
xmajorgrids,
yminorgrids,
ymajorgrids,
axis background/.style={fill=white},
legend style={at={(0.32,0.01)}, anchor=south west, legend cell align=left, align=left, draw=white!15!black, font=\footnotesize},
]
\addplot [color=brown, line width=1pt, mark size = 3pt,mark=square, mark options={solid, brown}]
  table[row sep=crcr]{%
9.6082e-06 	  1.7184e+02\\
3.4536e-06 	  2.1520e+02\\
1.5683e-06 	  2.5291e+02\\
7.4912e-07 	  2.9985e+02\\
4.3947e-07 	  3.4614e+02\\
};
\addlegendentry{\textsc{krogstad-v}}

\addplot [color=teal, line width=1pt, mark size = 3pt, mark=triangle, mark options={solid, teal}]
  table[row sep=crcr]{%
9.3515e-06 	 1.0544e+00\\
3.2531e-06 	 1.3415e+00\\
1.4198e-06 	 1.4868e+00\\
7.1699e-07 	 1.8227e+00\\
4.0083e-07 	 2.1668e+00\\
};
\addlegendentry{\textsc{krogstad-m}}

\end{axis}

\end{tikzpicture}%
}
    \caption{Results of Example~\ref{ex:source} in 2D using \textsc{krogstad-m}
    and \textsc{krogstad-v},
    see Section~\ref{sec:ex3_2D}. Left plot: increasing number of DOF per
    direction $n=[200,300,400,500,600]$, fixed number of time steps
    25. Right plot: increasing number of time steps $[15,20,25,30,35]$, fixed
    number of DOF per direction $n=400$.}
    \label{fig:ex3D2}
  \end{figure}

  We finally focus on the last two-dimensional experiment in MATLAB language
  and demonstrate that the proposed technique is effective also in the context
  of high-order methods. To this aim, we consider again
  Example~\ref{ex:source},
  but we now employ the fourth-order space discretization~\eqref{eq:FD4}
  of the fractional
  differential operator. For the time marching,
  we use the fourth-order Krogstad
  method in its matrix-oriented form
  (\textsc{krogstad-m}~\eqref{eq:krogstad-m}) and in its
  vector-oriented form (\textsc{krogstad-v}).
  Similarly to the previous experiment,
  the hyperparameters for shift-and-invert Lanczos are $m=10$
  and $\xi=\tau/10$,
  and \verb+pcg+ is run with default tolerance, maximum number of
  iterations 20, and initial guess the null vector.
  The results are collected in Figure~\ref{fig:ex3D2}.
  Here, the average number of PCG iterations
  for \textsc{krogstad-v} spans from $3.98$ to $4.70$ in the left plot, and
  from $4.71$ to $3.98$ in the right one.
  As expected, the proposed technique outperforms the state-of-the-art
  vector-oriented one
  also in this instance, with an even higher speedup
  compared to the previous examples  (roughly 170 in average).

\subsection{Three-dimensional experiments}\label{sec:numexp3}
In the previous section, we presented the two-dimensional experiments performed
in MATLAB language. Here, we focus our attention to some three-dimensional examples
on both consumer-level and professional hardware. In particular, our aim is
to show that the proposed technique scales favorably also in higher dimension
and on modern computer systems such as multi-core CPUs and GPUs.
This is indeed expected since, as already mentioned in Section~\ref{sec:ddim},
the realization of the needed operations in the introduced tensor framework
(i.e., the $\mu$-mode product and the Tucker operator) is strongly based on
BLAS (dense matrix-matrix products, to be precise), for which
heavily optimized routines exist nowadays.
We mention here, among the others, implementations provided by Intel MKL~\cite{mkl}
and OpenBLAS~\cite{xianyi2012model} for CPU-based systems and
cuBLAS~\cite{cublas} and MAGMA~\cite{ntd10_vecpar} for GPUs.
Also, (massive) parallelization can be directly exploited to evaluate
the nonlinear part of the system and the tensor $\boldsymbol F$ in
formula~\eqref{eq:fKu_diag}. In fact, in our context they are all independent
pointwise operations, for which GPUs are particularly well-suited.

We now summarize our tensor-oriented technique for realizing the linearized
BDF2 scheme~\eqref{eq:lbdf2-v} and the Strang splitting
method~\eqref{eq:Strang}.
For the former, we compute the first step as
\begin{subequations}
\begin{equation}
  \boldsymbol U_{1}=\left(\boldsymbol R_{\tau}\circ
  \left(\left(\boldsymbol U_0+
    \tau\boldsymbol G(0,\boldsymbol U_0)\right)\bigtimes_{\mu=1}^3
    Q_\mu^{\mathsf T}\right)\right)
    \bigtimes_{\mu=1}^3 Q_\mu
\end{equation}
and proceed for the following steps as
\begin{equation}\label{eq:lbdf2-t}
      \begin{aligned}
    \boldsymbol U_{n2}&=2\boldsymbol U_{n}-\boldsymbol U_{n-1},\\
    \boldsymbol U_{n+1}&=\left(
    \boldsymbol R_{2\tau/3}\circ \left(\left(    
    \tfrac{4}{3}\boldsymbol U_n
    -\tfrac{1}{3}\boldsymbol U_{n-1}
    +\tfrac{2\tau}{3}\boldsymbol G(t_{n+1},\boldsymbol U_{n2})\right)\bigtimes_{\mu=1}^3 Q_\mu^{\mathsf T}\right)\right)\bigtimes_{\mu=1}^3 Q_\mu.
  \end{aligned}
    \end{equation}
\end{subequations}
Here, we introduced the order-3 tensors
\begin{equation*}
  \boldsymbol R_\theta=\left((1-\theta\lambda_{j_1j_2j_3})^{-1}\right)_{j_1j_2j_3}\in
  \CC^{n\times n\times n},
  \quad \theta\in\left\{\tfrac{2\tau}{3},\tau\right\}.
\end{equation*}
This scheme will be referred to as \textsc{lbdf2-t}. Remark that the main
cost of the procedure is the computation of
two 3D Tucker operators
per time step.
The tensor-oriented Strang splitting method (denoted by \textsc{strang-t}) is
given by
\begin{equation}\label{eq:strang-t}
  \boldsymbol U_{n+1}=\Phi_{\tau/2}\left(\Phi_{\tau/2}(\boldsymbol U_n)
  \bigtimes_{\mu=1}^3E_\mu\right),
\end{equation}
where the small-sized matrix exponentials $E_\mu=\exp(\tau A_\mu)$ are computed once and for all before
the actual time integration starts and have in fact negligible computational
cost
(the main cost per time step is indeed the single 3D Tucker operator).

For the actual experiments on consumer-level hardware, in terms of software
and hardware we employ the ones already presented in Section~\ref{sec:numexp2}, i.e.,
MathWorks MATLAB\textsuperscript{\textregistered} on a standard laptop.
In addition, we implement the schemes in the C++
language, exploiting  the Intel OneAPI
MKL library (version 2021.4.0) for BLAS and OpenMP for basic parallelization.
The hardware is also equipped with a mobile NVIDIA card
(NVIDIA GeForce GTX 1650 Ti,
4GB of dedicated memory), which we exploit for the GPU experiments.
In this case, we use CUDA version 10.1, BLAS provided
by cuBLAS, and CUDA kernels for massive parallelization.
Even though this is clearly a consumer-level GPU, it can still be effectively
used for the relatively small-scale simulations being considered at this stage,
as we will demonstrate later on. Remark also that, in this setting, all the
experiments are performed in single-precision arithmetic for consistency with
the employed hardware.
On the other hand, for the experiments on professional hardware we use
double-precision arithmetic and exploit a workstation equipped with dual socket
Intel Xeon Gold 5118 with $2\times12$ cores as CPU, while the GPU is a
single NVIDIA V100 card (equipped with 16 GB of RAM).
When invoking BLAS on the CPU, we use the Intel OneAPI MKL library
(version 2020.1.0), whereas on the GPU we use cuBLAS from CUDA 11.2.
Parallelization is again handled with OpenMP on the CPU and CUDA kernels on the
GPU.

We finally remark that for the practical realization of the $\mu$-mode
products and the Tucker operators we use a ``loops-over-GEMMs'' approach.
This is particularly effective
in our HPC context since it avoids expensive permutations of large-sized tensors,
which have substantial computational burden due to the intrinsic alteration of the memory
layout (they are in fact memory-bound operations). To this aim, in
MathWorks MATLAB\textsuperscript{\textregistered}
we exploit the function \verb+pagemtimes+, while in C++ and CUDA we use the
relevant batched BLAS routines of Intel MKL and cuBLAS, respectively, from
the BLAS-like extensions. For more details, we invite an interested reader to
check the discussions in References~\cite{CCEOZ22,CCZ23bis,LBPSV15,C24}.
Notice also that the diagonalization is always performed on the CPU by means
of the
LAPACK routines \verb+{S,D}SYEV+, and then sent to the GPU. 
In fact, this computation is a small-sized task for which GPUs are not
well-suited. This is also the only point in which there is communication among
CPU and GPU: for the latter, all the remaining computations (the complete time
integration, in particular) are entirely performed on the device.

\subsubsection{First numerical experiment}\label{sec:ex1_3D}
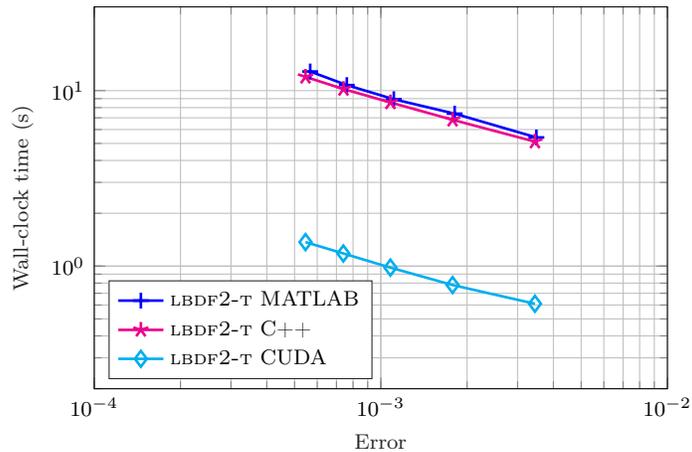
\begin{figure}[!htb]
  \centering
  \begin{tikzpicture}

\begin{axis}[%
width=3in,
height=2in,
font = \footnotesize,
scale only axis,
xmode=log,
xmin=1e-4,
xmax=1e-2,
xminorticks=true,
xlabel style={font=\color{white!15!black}},
xlabel = {\footnotesize Error},
xmajorgrids,
xminorgrids,
ymin=2e-1,
ymax=3e1,
ymode=log,
yminorticks=true,
ylabel style={font=\color{white!15!black}},
ylabel = {\footnotesize Wall-clock time (s)},
ymajorgrids,
yminorgrids,
axis background/.style={fill=white},
legend style={at={(0.025,0.025)}, anchor=south west, legend cell align=left, align=left, draw=white!15!black, font=\footnotesize},
]

\addplot [color=blue, line width=1pt, mark size = 3pt,mark=+, mark options={solid, blue}]
  table[row sep=crcr]{%
  3.49e-03  	5.41 \\
  1.81e-03  	7.39 \\
  1.11e-03  	8.93 \\
  7.61e-04  	10.75\\
  5.67e-04  	12.85\\
};
\addlegendentry{\textsc{lbdf2-t} MATLAB}

\addplot [color=magenta, line width=1pt, mark size = 3pt,mark=star, mark options={solid, magenta}]
  table[row sep=crcr]{%
  3.45e-03 	 5.12 \\
  1.78e-03 	 6.82 \\
  1.08e-03 	 8.56 \\
  7.41e-04 	 10.21\\
  5.46e-04 	 11.98\\
};
\addlegendentry{\textsc{lbdf2-t} C++}

\addplot [color=cyan, line width=1pt, mark size = 3pt,mark=diamond, mark options={solid, cyan}]
  table[row sep=crcr]{%
  3.45e-03 	 0.61\\
  1.78e-03 	 0.78\\
  1.08e-03 	 0.98\\
  7.40e-04 	 1.18\\
  5.46e-04 	 1.37\\
};
\addlegendentry{\textsc{lbdf2-t} CUDA}

\end{axis}

\end{tikzpicture}%
  \caption{Results of Example~\ref{ex:source} in 3D using \textsc{lbdf2-t}
  and increasing number of time steps $[15,20,25,30,35]$
  (fixed number of DOF per direction $n=200$),
  see Section~\ref{sec:ex1_3D}. The software architectures have been used with
  single-precision arithmetic and the hardware is a standard laptop.}
  \label{fig:cpugpu_ex1_3D}
\end{figure}

We start with Example~\ref{ex:source} and employ as time integrator
\textsc{lbdf2-t}~\eqref{eq:lbdf2-t}. The results of the experiment on the laptop
in single-precision arithmetic, for increasing number of time steps,
are summarized in Figure~\ref{fig:cpugpu_ex1_3D}.
As we can observe, the proposed technique has essentially the same performance
in MATLAB and in C++: this is expected, since both the architectures exploit
in a similar way the efficiency of BLAS and multithreading. On the other hand,
the advantage of using a GPU is clear. In fact, to obtain the same accuracy,
we have a computational speedup of a factor about 9.

\begin{figure}[!htb]
  \centering
  \begin{tikzpicture}

\begin{axis}[%
width=3in,
height=2in,
scale only axis,
font = \footnotesize,
xmode=log,
xmin=100,
xmax=800,
xminorticks=true,
xlabel style={font=\color{white!15!black}},
xmajorgrids,
xminorgrids,
ymajorgrids,
yminorgrids,
xtick={100,275,450,625,800},
xticklabels={100,275,450,625,800},
ymode=log,
ymin=1e-2,
ymax=7e2,
yminorticks=true,
ylabel style={font=\color{white!15!black}},
xlabel = {\footnotesize $n$},
ylabel = {\footnotesize Wall-clock time (s)},
axis background/.style={fill=white},
legend style={at={(0.015,0.8)}, anchor=south west, legend cell align=left, align=left, draw=white!15!black, font=\footnotesize},
]

\addplot [color=magenta, line width=1pt, mark size = 3pt, mark=star, mark options={solid, magenta}]
  table[row sep=crcr]{%
 200 	2.01e1\\
 275 	5.61e1\\
 350 	1.33e2\\
 425 	2.63e2\\
};
\addlegendentry{\textsc{lbdf2-t} C++}

\addplot [color=cyan, line width=1pt, mark size = 3pt, mark=diamond, mark options={solid, cyan}]
  table[row sep=crcr]{%
 200 	8.79e-1\\
 275 	2.55e0\\
 350 	6.09e0\\
 425 	1.32e1\\
};
\addlegendentry{\textsc{lbdf2-t} CUDA}

\addplot [color=magenta, dashed, line width=1pt, mark size = 3pt, mark=star, mark options={solid, magenta}]
  table[row sep=crcr]{%
 200 	3.63e-2\\
 275 	6.31e-2\\
 350 	9.05e-2\\
 425 	1.43e-1\\
};

\addplot [color=cyan, dashed, line width=1pt, mark size = 3pt, mark=diamond, mark options={solid, cyan}]
  table[row sep=crcr]{%
 200 	3.99e-2\\
 275 	6.41e-2\\
 350 	1.22e-1\\
 425 	1.48e-1\\
};

\end{axis}

\end{tikzpicture}%
  \caption{Results of Example~\ref{ex:source} in 3D using \textsc{lbdf2-t}
  and increasing number of DOF per direction $n=[200,275,350,425]$
  (fixed number of time steps 50),
  see Section~\ref{sec:ex1_3D}. The software architectures have been used with
  double-precision arithmetic and the hardware is a workstation. The dashed
  line represents the computational time needed to perform the diagonalization
  plus the communication time from CPU to GPU, when relevant.
  }
  \label{fig:cpugpu_ex1_3D_idof}
\end{figure}
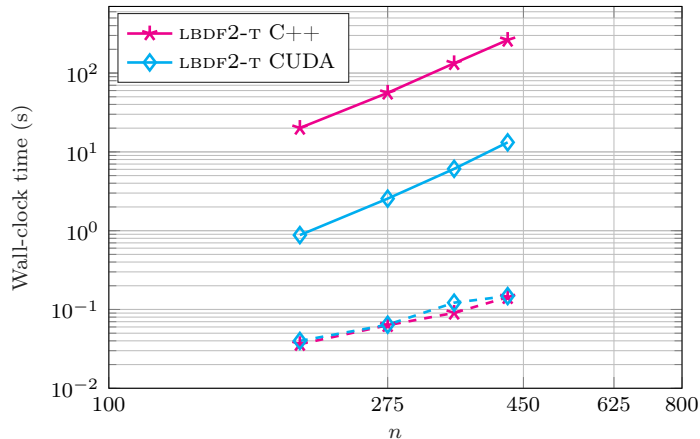

We then present in Figure~\ref{fig:cpugpu_ex1_3D_idof} the results of the experiment
with increasing number of DOF.
The time integrator is again \textsc{lbdf2-t}. For storage reasons,
and to show the performances of the scheme on professional hardware,
we perform the experiment
on the workstation in double-precision arithmetic.
The outcome is in line with the expectations, that is the GPU performance is
superior to the CPU one of more than an order of magnitude. Moreover, as mentioned
above, we clearly see that the time needed to compute the diagonalization
plus, when relevant, the sending of the result to the GPU is negligible
compared to the overall wall-clock time of the integration.

\subsubsection{Second numerical experiment}\label{sec:ex2_3D}
We now conclude with Example~\ref{ex:nosource} integrated in time
with \textsc{strang-t}~\eqref{eq:strang-t}.
The results of the experiment, for increasing number of time steps
in single-precision arithmetic on the laptop and for increasing number of DOF
in double-precision arithmetic
on the workstation, are summarized in Figures~\ref{fig:cpugpu_ex2_3D}
and~\ref{fig:cpugpu_ex2_3D_idof}, respectively.
\begin{figure}[!htb]
  \centering
  \begin{tikzpicture}

\definecolor{DarkGreen}{RGB}{1,90,32}

\begin{axis}[%
width=3in,
height=2in,
font=\footnotesize,
scale only axis,
xmode=log,
xmin=3e-4,
xmax=1e-1,
xminorticks=true,
xlabel style={font=\color{white!15!black}},
xlabel = {\footnotesize Error},
xmajorgrids,
xminorgrids,
ymajorgrids,
yminorgrids,
ymin=2e-1,
ymax=5e1,
ymode=log,
yminorticks=true,
ylabel style={font=\color{white!15!black}},
ylabel = {\footnotesize Wall-clock time (s)},
axis background/.style={fill=white},
legend style={at={(0.48,0.73)}, anchor=south west, legend cell align=left, align=left, draw=white!15!black, font=\footnotesize},
]

\addplot [color=purple, line width=1pt, mark size = 3pt, mark=asterisk, mark options={solid, purple}]
  table[row sep=crcr]{%
  1.57e-02 	 4.51 \\
  4.07e-03 	 7.93 \\
  1.79e-03 	 12.51\\
  9.77e-04 	 17.89\\
  6.03e-04 	 22.73\\
};
\addlegendentry{\textsc{strang-t} MATLAB}

\addplot [color=DarkGreen, line width=1pt, mark size = 3pt, mark=pentagon, mark options={solid, DarkGreen}]
  table[row sep=crcr]{%
  1.57e-02 	 4.22 \\
  4.12e-03 	 7.40 \\
  1.84e-03 	 11.75\\
  1.04e-03 	 16.42\\
  6.26e-04 	 20.88\\
};
\addlegendentry{\textsc{strang-t} C++}

\addplot [color=olive, line width=1pt, mark size = 3pt, mark=Mercedes star, mark options={solid, olive}]
  table[row sep=crcr]{%
  1.57e-02 	 0.45\\
  4.12e-03 	 0.76\\
  1.84e-03 	 1.14\\
  1.03e-03 	 1.52\\
  6.22e-04 	 1.91\\
};
\addlegendentry{\textsc{strang-t} CUDA}

\end{axis}

\end{tikzpicture}%
  \caption{Results of Example~\ref{ex:nosource} in 3D using \textsc{strang-t}
  and increasing number of time steps $[5,10,15,20,25]$
  (fixed number of DOF per direction $n=250$),
  see Section~\ref{sec:ex2_3D}. The software architectures have been used with
  single-precision arithmetic and the hardware is a standard laptop.}
  \label{fig:cpugpu_ex2_3D}
\end{figure}
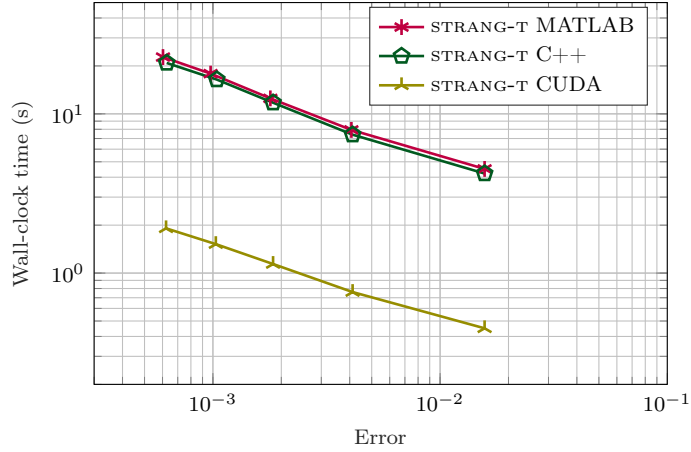%
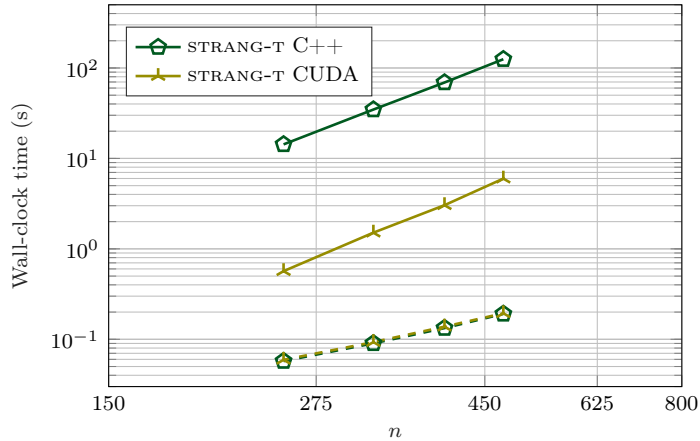
\begin{figure}[!htb]
  \centering
  \begin{tikzpicture}

\definecolor{DarkGreen}{RGB}{1,90,32}

\begin{axis}[%
width=3in,
height=2in,
font = \footnotesize,
scale only axis,
xmode=log,
xmin=150,
xmax=800,
xminorticks=true,
xlabel style={font=\color{white!15!black}},
xmajorgrids,
xminorgrids,
ymajorgrids,
yminorgrids,
xtick={150,275,450,625,800},
xticklabels={150,275,450,625,800},
ymode=log,
ymin=3e-2,
ymax=5e2,
yminorticks=true,
ylabel style={font=\color{white!15!black}},
xlabel = {\footnotesize $n$},
ylabel = {\footnotesize Wall-clock time (s)},
axis background/.style={fill=white},
legend style={at={(0.03,0.76)}, anchor=south west, legend cell align=left, align=left, draw=white!15!black, font=\footnotesize},
]

\addplot [color=DarkGreen, line width=1pt, mark size = 3pt, mark=pentagon, mark options={solid, DarkGreen}]
  table[row sep=crcr]{%
 250 	1.43e1\\
 325 	3.47e1\\
 400 	6.91e1\\
 475 	1.25e2\\
};
\addlegendentry{\textsc{strang-t} C++}

\addplot [color=olive, line width=1pt, mark size = 3pt, mark=Mercedes star, mark options={solid, olive}]
  table[row sep=crcr]{%
 250 	5.67e-1\\
 325 	1.51e0\\
 400 	3.04e0\\
 475 	5.98e0\\
};
\addlegendentry{\textsc{strang-t} CUDA}

\addplot [color=DarkGreen, dashed, line width=1pt, mark size = 3pt, mark=pentagon, mark options={solid, DarkGreen}]
  table[row sep=crcr]{%
 250 	5.74e-2\\
 325 	9.02e-2\\
 400 	1.33e-1\\
 475 	1.91e-1\\
};

\addplot [color=olive, dashed, line width=1pt, mark size = 3pt, mark=Mercedes star, mark options={solid, olive}]
  table[row sep=crcr]{%
 250 	5.93e-2\\
 325 	9.31e-2\\
 400 	1.38e-1\\
 475 	1.94e-1\\
};

\end{axis}

\end{tikzpicture}%
  \caption{Results of Example~\ref{ex:nosource} in 3D using \textsc{strang-t}
  and increasing number of DOF per direction $n=[250,325,400,475]$
  (fixed number of time steps 30),
  see Section~\ref{sec:ex2_3D}. The software architectures have been used with
  double-precision arithmetic and the hardware is a workstation. The dashed
  line represents the computational time needed to perform the diagonalization
  and computation of matrix exponentials
  plus the communication time from CPU to GPU, when relevant.
  }
  \label{fig:cpugpu_ex2_3D_idof}
\end{figure}%
The conclusions are essentially the same drawn for the previous experiment, that
is, the efficiency of GPUs compared to CPUs is neat in all instances.

\begin{figure}[!htb]
  \centering
  \includegraphics[scale=0.475]{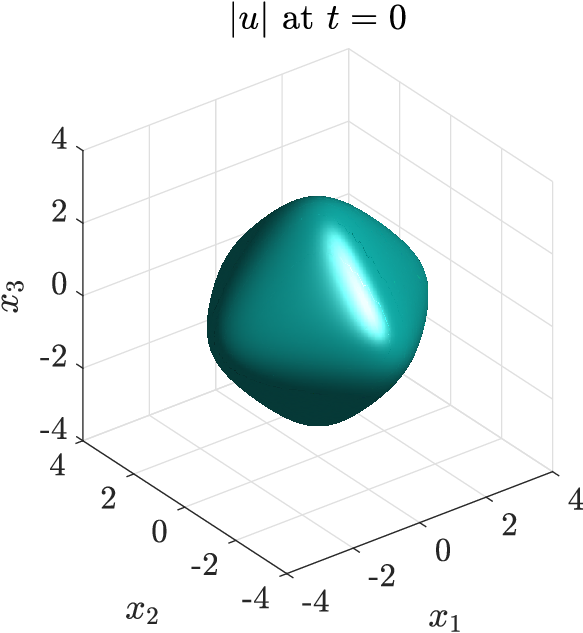}\hfill
  \includegraphics[scale=0.475]{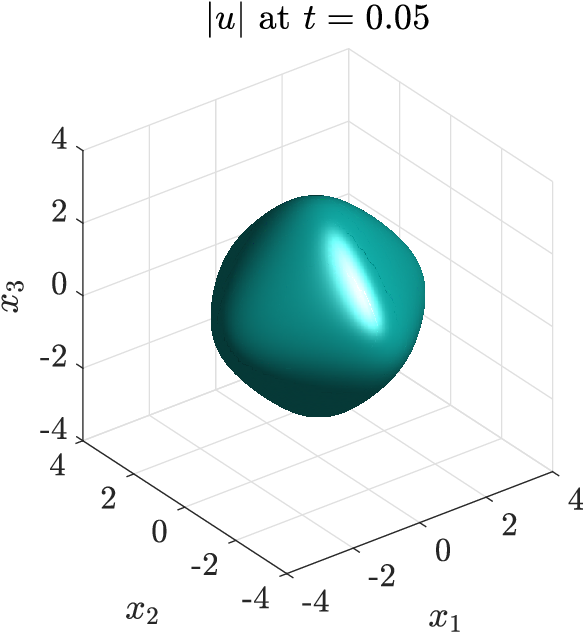}\\[2ex]
  \includegraphics[scale=0.475]{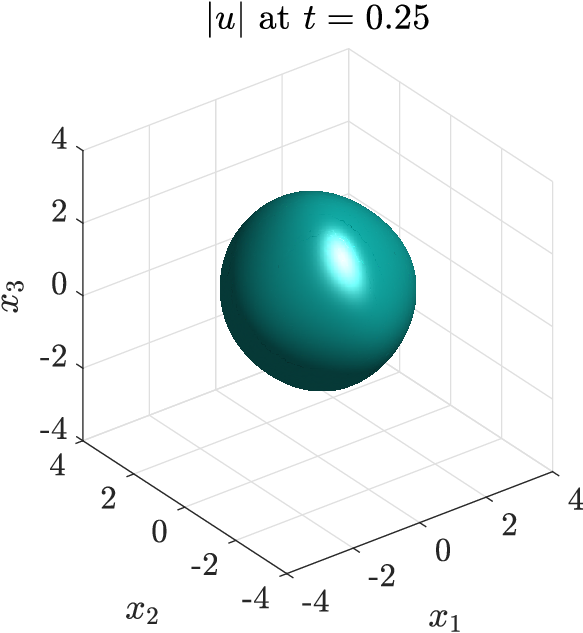}\hfill
  \includegraphics[scale=0.475]{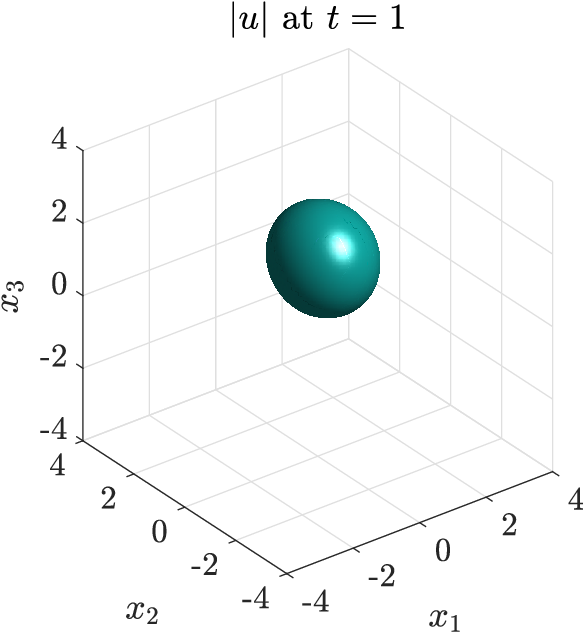}
  \caption{Isosurface of level $0.1$ of $\lvert u\rvert$ at different times
  for Example~\ref{ex:nosource}
  in 3D with time integrator \textsc{strang-t}, see Section~\ref{sec:ex2_3D}.}
  \label{fig:contour_3D}
\end{figure}

In addition,
we display the evolution of an isosurface of $\lvert u\rvert$ in
Figure~\ref{fig:contour_3D}, which shrinks and collapses as time proceeds.

\section{Conclusions and further developments}\label{sec:conc}
In this paper, we showed how to suitably employ matrix- and tensor-oriented
approaches to efficiently compute the numerical solution of the
multidimensional
evolutionary space-fractional complex Ginzburg--Landau equation.
After the semidiscretization of the spatial operator by appropriate finite
differences, the technique combines the arising
Kronecker sum structure and a diagonalization procedure
to compute in a fast way the actions of matrix functions needed by relevant
time-marching methods (e.g., linearly implicit and of exponential type).
We numerically demonstrated the superiority of the approach compared to
state-of-the-art techniques in a wide variety of two- and three-dimensional
situations, showing among the other things that multicore CPUs and GPUs can
be effectively exploited to accelerate the computations.

As interesting future works we plan to focus on other space-fractional
differential equations and develop an efficient, reliable, and stable technique
to be used when the arising discretization matrices are not diagonalizable
and/or not well-conditioned.

\section*{Acknowledgments}
The authors are members of the Gruppo Nazionale Calcolo
Scientifico-I\-sti\-tu\-to Na\-zio\-na\-le di Alta Matematica (GNCS-INdAM).
Fabio Cassini holds a post-doc fellowship funded by INdAM.
Also, the authors are grateful to Lukas Einkemmer and the Research Area Scientific
Computing of the University of Innsbruck for the availability to use their
professional hardware resources.
\bibliographystyle{elsarticle-harv}
  \bibliography{CC25}
\end{document}